\documentclass{article}
\usepackage{cite}
\usepackage{amsmath,amssymb,amsfonts}
\usepackage{algorithmic}
\usepackage[ruled,linesnumbered]{algorithm2e}
\usepackage{graphicx}
\usepackage{textcomp}
\usepackage{soul}
\usepackage{color}
\usepackage{colortbl}
\usepackage{lineno}
\usepackage{enumerate}
\usepackage{ntheorem}
\allowdisplaybreaks[4]
\newtheorem{theorem}{Theorem}{}
\newtheorem{corollary}{Corollary}{}
\newtheorem{remark}{Remark}{}
\newtheorem{lemma}{Lemma}{}
\newtheorem{assumption}{Assumption}{}
\newtheorem{definition}{Definition}{}
\newtheorem*{proof}{Proof}{}
\title{Accelerated Distributed Aggregative Optimization}
\author{
Jiaxu Liu$^a$, Song Chen$^a$, Shengze Cai$^b$,  Chao Xu$^b$, Jian Chu$^b$\thanks{Corresponding author: shengze\_cai@zju.edu.cn.} \\
~\\
$^a$ School of Mathematical Sciences, Zhejiang University, Hangzhou, 310027, China  \\
$^b$ Institute of Cyber-Systems \& Control, 
College of Control Science \& Engineering, \\
Zhejiang University, Hangzhou, 310027, China  \\
}

\begin{document}

\maketitle
\begin{abstract}
This paper delves into the investigation of a distributed aggregative optimization problem within a network. In this scenario, each agent possesses its own local cost function, which relies not only on the local state variable but also on an aggregated function of state variables from all agents. To expedite the optimization process, we amalgamate the heavy ball and Nesterov’s accelerated method with distributed aggregative gradient tracking, resulting in the proposal of two innovative algorithms, aimed at resolving the distributed aggregative optimization problem. Our analysis demonstrates that the proposed algorithms can converge to an optimal solution at a global linear convergence rate when the objective function is  strongly convex with the Lipschitz-continuous gradient, and when the parameters (e.g., step size and momentum coefficients) are chosen within specific ranges. Additionally, we present several numerical experiments  to verify the effectiveness, robustness and superiority of our proposed algorithms.

\end{abstract}

\section{Introduction}
Distributed optimization plays a crucial role in machine learning, particularly in scenarios where the dataset or model is too large to fit in a single machine. By distributing the dataset or model across parallel machines and employing distributed optimization techniques, one can achieve more efficient computation and scalable machine learning \cite{dean2012large,Barbarossa2014heterogeneous,li2014communication,Predd2009ACT}. Furthermore, distributed optimization has garnered significant attention in the field of control due to its broad applications, including formation control \cite{WANG2014nonholonomic}, sensor networks \cite{Zhu2013Sensor}, resource allocation \cite{Deng2018DistributedCA}, and more. In a distributed optimization problem, each agent can only access its own information and communicate data with its neighbors to cooperatively minimize a global cost function. Through communication and collaboration, these agents can work in parallel, resulting in a faster optimization process that can handle larger datasets{\cite{lee2014model}}.

To address the distributed optimization problem, several remarkable algorithms have been developed. Many well-known distributed optimization algorithms, based on a consensus scheme, employ gradient descent techniques. These include the distributed subgradient descent \cite{Nedic2009Subgradient}, the distributed dual averaging gradient algorithm \cite{Duchi2010DualAF}, and the push-sum distributed algorithm \cite{Tsianos2012PushSumDD}, among others. However, these algorithms suffer from slow convergence due to the use of a gradually vanishing step size in their design. Although employing a constant step size can improve the convergence rate of distributed gradient-based algorithms, they are limited to converging only to a small neighborhood of the optimal solution due to the use of local gradients in each agent \cite{yuan2016convergence}. To address the issues of slow convergence rate and non-optimal solutions, the gradient tracking strategy has been integrated into distributed convex optimization algorithms. These algorithms utilized an estimate of the global average gradient to replace the local gradient in each agent \cite{xu2015augmented,Nedi2016AchievingGC}. In \cite{shi2015extra}, a  decentralized exact first-order algorithm (EXTRA) was proposed, which achieved geometric convergence to the global optimal solution by introducing a cumulative correction term.

Moreover, to expedite the convergence rate, acceleration algorithms such as the heavy ball method \cite{polyak1964some} and Nesterov’s accelerated method \cite{nesterov1983method} can be employed in distributed optimization. In \cite{jakovetic2014fast}, the distributed Nesterov’s gradient method (D-NG) was proposed, improving the convergence rate to $O(\frac{\log k}{k})$. Additionally, in \cite{xin2019distributed}, two distributed Nesterov's gradient methods over strongly connected and directed networks were introduced by extending $\mathcal{A}\mathcal{B}$  \cite{Ran2018Geometric} with Nesterov’s momentum. Subsequently, in \cite{Ran2020Generalization}, the authors combined $\mathcal{A}\mathcal{B}$ 
 with the heavy ball method, proposing a distributed heavy-ball algorithm named $\mathcal{A}\mathcal{B}m$ . It was demonstrated that $\mathcal{A}\mathcal{B}m$ 
 exhibits a global 
linear rate when the step-size and momentum parameters were positive and sufficiently small. Furthermore, Qu and Li \cite{Qu2020Accelerated} further integrated gradient tracking with distributed Nesterov’s gradient descent method, resulting in two accelerated distributed Nesterov’s methods termed as Acc-DNGD-SC and Acc-DNGD-NSC. Additionally, a distributed Nesterov-like gradient tracking algorithm (D-DNGT) was introduced in \cite{lu2020nesterov}, which incorporated gradient tracking into the distributed Nesterov’s method with momentum terms and employed nonuniform step sizes. These algorithms achieved linear convergence rates for smooth and strongly convex objective functions, significantly enhancing the convergence speed compared to algorithms without momentum terms.

In the distributed optimization problems mentioned above, a prevalent setup that has emerged is known as consensus optimization. However, in various other scenarios, such as the multi-agent formation control problem, multi-robot surveillance scenario \cite{Carnevale2022AggregativeFO}, and aggregative game \cite{koshal2016distributed}, the objective function of each agent is not only dependent on its local state but also determined by other agents’ variables through an aggregative variable. Such an optimization problem is referred to as the distributed aggregative optimization in \cite{Li2022Aggregative}, where a distributed aggregative gradient tracking (DAGT) method was proposed to address this challenge. Following \cite{Li2022Aggregative}, Li et al. \cite{Li2022Online} considered online convex optimization with an aggregative variable to solve time-varying cost functions. Subsequently, the online distributed aggregative optimization problem with constraints was investigated in \cite{Carnevale2022Coordination}. Moreover, Chen and Liang \cite{Chen2022DistributedAO} considered finite bits communication and proposed a novel distributed aggregative optimization algorithm with quantification scheme. In \cite{Wang2022DistributedPA}, the distributed convex aggregative optimization was further combined with the Frank-Wolfe algorithm to solve the aggregative optimization problem over time-varying communication graphs.

To the best of our knowledge,  little work has been
done to accelerate the convergence rate of distributed aggregative optimization. Therefore, we propose two novel accelerated algorithms for distributed aggregative optimization problem in this paper. The main contributions of our work can be summarized as follows. 
\begin{enumerate}
    \item  We combine accelerated algorithms with distributed aggregative gradient tracking (DAGT)\cite{Li2022Aggregative}, resulting in distributed aggregative gradient tracking with the heavy ball algorithm (DAGT-HB) and distributed aggregative gradient tracking with Nesterov's algorithm (DAGT-NES).
    \item  We provide theoretical analysis demonstrating that the DAGT-HB and DAGT-NES algorithms can converge to an optimal solution at a global linear convergence rate { under several standard assumptions.} Moreover, we offer guidelines for selecting suitable parameter ranges, including the step size and momentum term. Furthermore, within the context of the quadratic distributed aggregative optimization problem, we illustrate that the convergence rates of DAGT-HB and DAGT-NES surpass that of DAGT when the spectral radius of the weighted adjacency matrix falls within a specific range.
    \item  We perform numerical simulations to verify the effectiveness, robustness and superiority of the proposed DAGT-HB and DAGT-NES algorithms, which further support our theoretical findings presented in this paper.

\end{enumerate}

The remainder of this paper is structured as follows. In Section \uppercase\expandafter{\romannumeral2},  the  basic notations, the basic definition of graph and the distributed aggregative optimization problem are presented. Sections \uppercase\expandafter{\romannumeral3} and \uppercase\expandafter{\romannumeral4} are dedicated to the introduction of DAGT-HB and DAGT-NES, along with the analysis of their convergence rates. Section \uppercase\expandafter{\romannumeral5} delves into the analysis and comparison of convergence rates among DAGT, DAGT-HB, and DAGT-NES for the quadratic distributed aggregative optimization problem. To validate the proposed algorithms, several numerical examples are provided in Section \uppercase\expandafter{\romannumeral6}. Finally, Section \uppercase\expandafter{\romannumeral7} serves as the conclusion of the paper.

\vspace{-1 em}
\section{Preliminaries And Problem Setup}
\subsection{Basic Notations}
The set of real and positive real numbers are denoted by $\mathbb{R}$ and $\mathbb{R}^{\mathbb{+}}$, the set of $n$-dimensional column vectors is denoted by $\mathbb{R}^n$. Let $\mathbf{1}_n$ and $\mathbf{0}_n$ represent the column vectors of $n$ ones and zeros, respectively. $I_n$ denotes the $n\times n$  dimensional identity matrix. Let $||\cdot||$ and $x^\top$ be the standard Euclidean norm (or induced matrix norm) and the transpose of $x\in \mathbb{R}^n $.
We use pointwise order $>$ for any vectors $a=(a_1, a_2, \ldots, a_n)^\top$ and $b=(b_1, b_2, \ldots, b_n)^\top \in \mathbb{R}^n$, i.e., $a> b \iff a_i>b_i$,  $i=1, 2, \ldots, n$. For vectors $x_1, \cdots, x_N \in \mathbb{R}^n$, we use the notation $x=\operatorname{col}\left(x_1, \cdots\right.$, $\left.x_N\right)=\left(x_1^{\top}, \cdots, x_N^{\top} \right)^{\top}$ to denote a new stacked vector. Also, we  use
blkdiag$(A_1, A_2, \ldots, A_N)$ to represent the block diagonal matrix
where the $i$-th diagonal block is given by the matrix $A_i \in \mathbb{R}^{m_i\times n_i}$, $i=1, 2, \ldots, N$. The Kronecker product of arbitrary matrices $A \in \mathbb{R}^{m \times n}$ and $B \in \mathbb{R}^{p \times q}$ is defined as $A \otimes B \in$ $\mathbb{R}^{m p \times n q}$. Let $\rho(A)$ be the spectral radius of a square matrix $A$. $A>0$ represents that $A$ is positive, that is, every entry of $A$ is greater than 0. Let $K_N=\frac{1}{N}\mathbf{1}_N\mathbf{1}^\top_N$. Besides, $\nabla f(\cdot)$ is the gradient of a differentiable function $f(\cdot)$.
\vspace{-1 em}
\subsection{Graph Theory}
Here, we provide some basic definitions of graph theory. {Let $\mathcal{G}=(\mathcal{V}, \mathcal{E}, A)$ denote a weighted undirected graph with a finite vertex set $\mathcal{V}=\{1, \ldots, N\}$, an edge set $\mathcal{E} \subseteq \mathcal{V} \times \mathcal{V}$, and a weighted adjacency matrix $A=\left[a_{i j}\right] \in \mathbb{R}^{N \times N}$ with $a_{i j}=a_{ji}>0$ if $(j, i) \in \mathcal{E}$ and $a_{i j}=a_{ji}=0$ otherwise.} $\mathcal{N}_i=\{j:(j, i) \in \mathcal{E}\}$ denotes the set of neighbors of agent $i$ and $d_i=\sum_{j=1}^N a_{i j}$ denotes the weighted degree of vertex $i$. The
graph $\mathcal{G}$ is called connected if for any $i, j \in  \mathcal{V}$,  there exists a  path from
$i$ to $j$. The Laplacian matrix of graph $\mathcal{G}$ is $L=D-A$ with $D=\operatorname{diag}\left\{d_1, \ldots, d_N\right\}$. The real eigenvalues of $L$ are denoted by $\lambda_1,...,\lambda_N$ with $\lambda_i \leq \lambda_{i+1}$, $i=1, \ldots, N-1$. The undirected graph $\mathcal{G}$ is connected if and only if $\lambda_2>0$.
\vspace{-1 em}
\subsection{Problem Formulations}
{In this paper, we consider the distributed aggregative
optimization problem that can be written as:
\begin{align}
    \min _{x \in \mathbb{R}^n} F(x), F(x)  =\sum_{i=1}^N f_i\left(x_i, u(x)\right) ,u(x)  =\frac{\sum_{i=1}^N \phi_i\left(x_i\right)}{N},
\end{align}
where $x=\operatorname{col}\left(x_1, \ldots, x_N\right)$ is the global state variable with $x_i \in$ $\mathbb{R}^{n_i}, n=\sum_{i=1}^N n_i$, and {$f_i: \mathbb{R}^{n_i}\times \mathbb{R}^{d} \rightarrow \mathbb{R}$} is the local objective function.} In problem (1), $u(x) \in \mathbb{R}^{d}$  represents an aggregative variable that can access information from all agents, while the function $\phi_i$: $\mathbb{R}^{n_i} \rightarrow \mathbb{R}^d$ is only accessible to agent $i$. Additionally, each agent $i$
 has knowledge of its own state variable  $x_i$
 but can not access information about other agents’ state variables. Moreover,  each agent $i$ can only privately access the information on $f_i$. The purpose of this paper is to develop a distributed optimization algorithm to {solve the problem (1).}


 For simplicity in the following analysis, we  denote  $\nabla_1 f_{i}(\cdot, \cdot)$ and  $\nabla_2 f_{i}(\cdot, \cdot)$ as the gradient of $f_{i}$ with respect to the first argument and the second argument respectively  for all $i\in \{1,2,\cdots,N\}$.  For $x=$ col$\left(x_1, \ldots, x_N\right) \in \mathbb{R}^{n}$ and $y=$ col$\left(y_1, \ldots, y_N\right) \in \mathbb{R}^{N d}$, we define $f(x, y)=\sum_{i=1}^N f_i\left(x_i, y_i\right)$, $\nabla_1 f(x, y)=$ col$\left(\nabla_1 f_1\left(x_1, y_1\right), \ldots, \nabla_1 f_N\left(x_N, y_N\right)\right)$ and $\nabla_2 f(x, y)=$ col$\left(\nabla_2 f_1\left(x_1,y_1\right), \ldots, \nabla_2 f_N\left(x_N, y_N\right)\right)$. Also, we denote $\phi(x)=$ col$(\phi_1(x_1),\ldots,\phi_N(x_N))$ and   $\nabla \phi(x)=$ blkdiag$(\nabla\phi_1\left(x_1\right), \ldots, \nabla\phi_N\left(x_N\right))\in \mathbb{R}^{n\times Nd}$. Hence, the gradient of $F(x)$ is described by
{$$
\nabla F(x)=\nabla_1 f\left(x, \mathbf{1}_N \otimes u(x)\right)+\nabla \phi(x) \mathbf{1}_N \otimes \frac{1}{N} \sum_{i=1}^{{N}} \nabla_2 f_i\left(x_i, u(x)\right).
$$}



To facilitate the subsequent analysis, it is necessary to make some definitions and general assumptions. Firstly, we give the definitions of $L$-smooth function and $\mu$-strongly convex function respectively.

\begin{definition}
   {\cite{nesterov2003introductory}} A differentiable function $f: \mathbb{R}^n \rightarrow \mathbb{R}$ is $L$-smooth if for all $x, y \in \mathbb{R}^n$
\begin{equation}
    \|\nabla f(x)-\nabla f(y)\| \leq L\|x-y\|.
\end{equation}
Following \cite{nesterov2003introductory}, it is equivalent to that for all $x, y \in \mathbb{R}^n$
\begin{equation}
    f(y) \leq f(x)+\nabla f(x)^{\top}(y-x)+\frac{L}{2}\|y-x\|^2.
\end{equation}
\end{definition}

\begin{definition}
    \cite{nesterov2003introductory} A differentiable function $f: \mathbb{R}^n \rightarrow \mathbb{R}$ is $\mu$-strongly convex  if for all $x, y \in \mathbb{R}^n$
\begin{equation}
    \quad \mu\|x-y\|^2 \leq(x-y)^{\top}(\nabla f(x)-\nabla f(y)). 
\end{equation}
 It also implies that for all $x, y \in \mathbb{R}^n$
\begin{equation}
    f(x)+\nabla f(x)^{\top}(y-x)+\frac{\mu}{2}\|y-x\|^2 \leq f(y).  
\end{equation}
\end{definition}


Next, we make some common assumptions about the graph $\mathcal{G}$  and objective function $F(x)$.

\begin{assumption}
    The undirected graph $\mathcal{G}$ is connected and $A$ is  symmetric and doubly stochastic, that is, $\sum_{i=1}^{N}a_{ij}=1$ and  $\sum_{j=1}^{N}a_{ij}=1$ for all $i,j=1,\ldots,N$.
\end{assumption}

{\begin{assumption}
    The global objective function $F$ is  $\mu$-strongly convex. For all $i\in \mathcal{V}$, $\nabla_1 f_i(x_i, y_i)+\nabla \phi_i(x_i) \frac{1}{N} \sum_{i=1}^N \nabla_2 f_i\left(x_i, y_i\right)$ is  $L_{1i}$-Lipschitz continuous. Namely, for all $x,x'  \in \mathbb{R}^n$ and $y,y'  \in \mathbb{R}^{Nd}$,
$$
    \begin{aligned}
        &\left\|\nabla_1 f(x, y)+\nabla \phi(x) \mathbf{1}_N \otimes \frac{1}{N} \sum_{i=1}^N \nabla_2 f_i\left(x_i, y_i\right)-\nabla_1 f(x', y')\right.\\
        &\left.-\nabla \phi(x') \mathbf{1}_N \otimes \frac{1}{N} \sum_{i=1}^N \nabla_2 f_i\left(x_i', y_i' \right)  \right\|\\
        & \leq L_1 \left( \|x-x'\|+\|y-y'\|\right),
    \end{aligned}
$$
where $L_1=\max\{L_{11},\ldots,L_{1N}\}$, which can deduce that the global objective function $F$ is $L_1$-smooth.
\end{assumption}

\begin{assumption}
  For all $i\in \mathcal{V}$,   $\nabla_2 f_i(x_i, y_i)$ is $L_{2i}$-Lipschitz continuous. Namely, for all $x,x'  \in \mathbb{R}^n$ and $y,y'  \in \mathbb{R}^{Nd}$, $\|\nabla_2 f(x,y)-\nabla_2 f(x',y')\|  \leq L_2 \left( \|x-x'\|+\|y-y'\|\right)$, where $L_2=\max\{L_{21},\ldots,L_{2N}\}$.
\end{assumption}

\begin{assumption}
    For all $i=1,2,\ldots,N$, the aggregation function $\phi_i$ is differentiable and 
    $\|\nabla \phi_i(x)\|$ is bounded by  $L_{3i}>0$, which implies $\|\nabla \phi(x)\|\leq L_3$ and $\phi$  is $L_3$-Lipschitz continuous, where  $L_3=\max\{L_{31},\ldots,L_{3N}\}$.
\end{assumption}
}

\begin{remark}
Assumptions 2-4 are general and the same as \cite{Li2022Aggregative}\cite{Chen2022DistributedAO} \cite{wang2024momentum}.
\end{remark}

To perform  the following analysis, several key lemmas are listed below.

{
\begin{lemma}
\cite{Li2022Aggregative} Let $f: \mathbb{R}^n \rightarrow \mathbb{R}$ be $\mu$-strongly convex and $L$-smooth.  When  $\alpha \in(0,1 / L]$, then for all $x, y \in \mathbb{R}^n$, we can obtain $\|x-\alpha \nabla f(x)-(y-\alpha \nabla f(y))\| \leq(1-\mu \alpha)\|x-y\|$.
\label{lem1}
\end{lemma}
}

\begin{lemma}
    
  {\cite{xiao2004fast} }Under Assumption 1,  for the adjacency
matrix $A$, the following properties hold:
\begin{enumerate}
    \item $\mathcal{A} \mathcal{K}=\mathcal{K} \mathcal{A}=\mathcal{K}$, where $\mathcal{A}=A \otimes I_d$, $\mathcal{K}=K_N \otimes I_d$.
    \item $\|\mathcal{A} x-\mathcal{K} x\| \leq \rho\|x-\mathcal{K} x\|$ for any $x \in$ $\mathbb{R}^{N d}$ and $\rho =\|A-K_N\|<1$.
    \item $||A-I_N|| \leq 2$.
\end{enumerate}
\label{lem2}
\end{lemma}

\begin{lemma}
    
 {(\cite{horn2012matrix} Corollary 8.1.29)} Let $W\in \mathbb{R}^{n\times n} $ be nonnegative and $x\in \mathbb{R}^{n}$ be positive. If $Wx<\lambda x$ with $\lambda>0$, then $\rho(W)<\lambda$.
 \label{lem3}
\end{lemma}

{\begin{lemma}
{(\cite{horn2012matrix} Lemma 5.6.10)} Let $W\in \mathbb{R}^{n\times n}$ and $\epsilon>0$ be given. There is a matrix norm $\| |\cdot|\|$ such
that $\rho(W)\leq \| |W|\|\leq \rho(W)+\epsilon$.
\label{lem4}
\end{lemma}}

\begin{lemma}
\cite{zheng2010extended} Let $\mathbb{R}$ be the real number field, and $H\left(z\right)$ denote the degree $n$ real coefficient polynomial
$$
H\left(z\right)=a_n z^n+a_{n-1} z^{n-1}+\cdots+a_2 z^2+a_1 z+a_0,
$$
 where $a_0, a_1, \ldots, a_n \in \mathbb{R}$.
The Jury matrix of $H\left(z\right)$ can be written as

\begin{tabular}{c|lllllll} 
& $z^0$ & $z^1$ & $z^2$ & $\cdots$ & $z^{n-2}$ & $z^{n-1}$ & $z^n$ \\
\hline 1 & $a_0$ & $a_1$ & $a_2$ & $\cdots$  & $a_{n-2}$ & $a_{n-1}$ & $a_n$ \\
2 & $a_n$ & $a_{n-1}$ & $a_{n-2}$ & $\cdots$ & $a_2$ & $a_1$ & $a_0$ \\
3 & $b_0$ & $b_1$ & $b_2$ & $\cdots$  & $b_{n-2}$ & $b_{n-1}$ & \\
4 & $b_{n-1}$ & $b_{n-2}$ & $b_{n-3}$ & $\cdots$  & $b_1$ & $b_0$ & \\
5 & $c_0$ & $c_1$ & $c_2$ & $\cdots$  & $c_{n-2}$ & \\
6 & $c_{n-2}$ & $c_{n-3}$ & $c_{n-4}$ & $\cdots$  & $c_0$ & \\
$\vdots$ & $\vdots$ & $\vdots$ & $\vdots$ & & & &  \\
$2 n-5$ & $l_0$ & $l_1$ & $l_2$ & $l_3$ & & &  \\
$2 n-4$ & $l_3$ & $l_2$ & $l_1$ & $l_0$ & & &  \\
$2 n-3$ & $m_0$ & $m_1$ & $m_2$ & & & &  \\
\end{tabular}

where
$$
\begin{aligned}
& b_i=\left|\begin{array}{ll}
a_0 & a_{n-i} \\
a_n & a_i
\end{array}\right|, \quad i=0,1,2, \ldots, n-1,\\ & c_j=\left|\begin{array}{ll}
b_0 & b_{n-j-1} \\
b_{n-1} & b_j
\end{array}\right|, \quad j=0,1,2, \ldots, n-2, \\
 & \cdots \\&  m_0=\left|\begin{array}{ll}
l_0 & l_3 \\
l_3 & l_0
\end{array}\right|, \quad m_1=\left|\begin{array}{cc}
l_0 & l_2 \\
l_3 & l_1
\end{array}\right|, \quad m_2=\left|\begin{array}{cc}
l_0 & l_1 \\
l_3 & l_2
\end{array}\right| . 
&
\end{aligned}
$$
All the {module} of roots of a real coefficient polynomial $H\left(z \right)\left(n \geqslant 3, a_n>\right.$ 0) are less than 1 if and only if the following four conditions hold:
\begin{enumerate}[1)]
    \item $H\left(1 \right)>0$;
    \item $(-1)^n H\left(-1 \right)>0$;
    \item $\left|a_0\right|<a_n$;
    \item $\left|b_0\right|>\left|b_{n-1}\right|, \left|c_0\right|>\left|c_{n-2}\right|, \ldots, \left|l_0\right|>\left|l_3\right|, \left|m_0\right|>\left|m_2\right|$.
\end{enumerate}
\label{lem5}
\end{lemma}

\vspace{-1 em}
\section{DAGT-HB}

In order to solve problem (1), we combine the distributed aggregative gradient tracking (DAGT)\cite{Li2022Aggregative} with heavy ball method and propose the following DAGT-HB algorithm:
\begin{align}
    x_{i, k+1} & =x_{i, k}-\alpha\left[\nabla_1 f_i\left(x_{i, k}, u_{i, k}\right)+\nabla \phi_i\left(x_{i, k}\right) s_{i, k}\right] \notag\\
& \quad +\beta (x_{i,k}-x_{i,k-1}), \\
u_{i, k+1} & =\sum_{j=1}^N a_{i j} u_{j, k}+\phi_i\left(x_{i, k+1}\right)-\phi_i\left(x_{i, k}\right), \\
s_{i, k+1} & =\sum_{j=1}^N a_{i j} s_{j, k}+\nabla_2 f_i\left(x_{i, k+1}, u_{i, k+1}\right)  \notag\\
&\quad  -\nabla_2 f_i\left(x_{i, k}, u_{i, k}\right),
\end{align}
{where $x_{i,k} \in \mathbb{R}^{n_i}$ and $u_{i,k},s_{i,k} \in \mathbb{R}^{d}$.}

In DAGT-HB, we incorporate a momentum term $x_{i,k}-x_{i,k-1}$
 to expedite the convergence rate of the algorithm. Given that distributed aggregative optimization has moved away from consistency protocols and interlinked agents through an aggregative variable, each agent only needs to reach the optimal point of its own local objective function. Hence, it is intuitive to consider that for each agent, if the iteration point generated by the algorithm consistently progresses towards the optimal solution, introducing momentum in the same direction can inevitably accelerate the convergence speed of the algorithm. As $u(x)$ 
 represents global information inaccessible to all agents, we introduce  $u_{i,k}$
 for agent $i$
 to track the average 
$u(x)$. Simultaneously, $s_{i,k}$  
 is introduced to track the gradient sum $\frac{1}{N}\sum_{i=1}^N\nabla_2 f_i\left(x_{i}, u(x)\right)$, which also can not be obtained by all agents.

To facilitate subsequent analysis, we define
$$
\begin{aligned}
    &x_k= \operatorname{col}\left(x_{1, k}, \ldots, x_{N, k}\right),u_k= \operatorname{col}\left(u_{1, k}, \ldots, u_{N, k}\right),\\
    &s_k= \operatorname{col}\left(s_{1, k}, \ldots, s_{N, k}\right), \phi(x_k)=\operatorname{col}\left(\phi_1(x_{1, k}),\ldots,\phi_N(x_{N, k})\right),\\
    &\nabla_1 f\left(x_k, u_k\right)=\operatorname{col}\left(\nabla_1 f_1\left(x_{1, k}, u_{1, k}\right),\ldots,\nabla_1 f_N\left(x_{N, k}, u_{N, k}\right)\right),\\
     &\nabla_2 f\left(x_k, u_k\right)=\operatorname{col}\left(\nabla_2 f_1\left(x_{1, k}, u_{1, k}\right),\ldots,\nabla_2 f_N\left(x_{N, k}, u_{N, k}\right)\right),\\
&\nabla\phi\left(x_k\right)=\operatorname{blkdiag}\left(\nabla\phi_1(x_{1, k}),\ldots,\nabla\phi_N(x_{N, k})\right).
\end{aligned}
$$
Then the DAGT-HB algorithm can be rewritten as the following compact form:
\begin{align}
    x_{k+1} & =x_k-\alpha\left[\nabla_1 f\left(x_k, u_k\right)+\nabla \phi\left(x_k\right) s_k\right]+\beta(x_k-x_{k-1}),\label{compact hx} \\
u_{k+1} & =\mathcal{A} u_k+\phi\left(x_{k+1}\right)-\phi\left(x_k\right), \label{compact hu}\\
s_{k+1} & =\mathcal{A} s_k+\nabla_2 f\left(x_{k+1}, u_{k+1}\right)-\nabla_2 f\left(x_k, u_k\right),\label{compact hs}
\end{align}
with $\mathcal{A}=A \otimes I_d \in \mathbb{R}^{Nd\times Nd}$ as defined in Lemma $2$. 


{Firstly we notice that multiplying $\frac{\mathbf{1}_{N}^{\top}}{N}\otimes I_{d}$ on both sides of \eqref{compact hu} and \eqref{compact hs} can lead to
\begin{align}
    &\bar{u}_{k+1}=\bar{u}_k+\frac{1}{N} \sum_{i=1}^N \phi_i\left(x_{i, k+1}\right)-\frac{1}{N} \sum_{i=1}^N \phi_i\left(x_{i, k}\right), \\
&\bar{s}_{k+1}=\bar{s}_k+\frac{1}{N} \sum_{i=1}^N\left[\nabla_2 f_i\left(x_{i,k+1}, u_{i,k+1}\right) -\nabla_2 f_i\left(x_{i,k}, u_{i,k}\right)\right],
\end{align}
where $\bar{u}_{k+1} \in \mathbb{R}^d$ and $\bar{s}_{k+1} \in \mathbb{R}^d$.} Then if  we initialize $u$ and $s$ as $u_{i,0}=\phi_i(x_{i,0})$ and $s_{i,0}=\nabla_2 f_i(x_{i,0},u_{i,0})$ for $i=1,2,\ldots,N$, where $x_{i,0}$ is arbitrary, we can derive 
\begin{align}
    \bar{u}_k & =\frac{1}{N} \sum_{i=1}^N u_{i, k}=\frac{1}{N} \sum_{i=1}^N \phi_i\left(x_{i, k}\right), \label{mean hu}\\
\bar{s}_k & =\frac{1}{N} \sum_{i=1}^N s_{i, k}=\frac{1}{N} \sum_{i=1}^N \nabla_2 f_i\left(x_{i, k}, u_{i, k}\right) \label{mean hs}.
\end{align}

Next, we establish the equivalence of the optimal solution to the problem (1) and the fixed point of the DAGT-HB algorithm. 

\begin{lemma}
    Under Assumption $1$ and Assumption $2$, the fixed point of \eqref{compact hx}-\eqref{compact hs} is
the optimal solution to the problem (1).
\end{lemma}

\begin{proof}Denote the equilibrium point of \eqref{compact hx}-\eqref{compact hs} as $x^*=$col $\left(x_1^*, \ldots, x_N^*\right)$, $u^*=$ col$\left(u_1^*, \ldots, u_N^*\right)$, and $s^*=\operatorname{col}\left(s_1^*, \ldots, s_N^*\right)$. From  \eqref{compact hx}-\eqref{compact hs} , we can obtain
\begin{align}
    & \nabla_1 f\left(x^*, u^*\right)+\nabla \phi\left(x^*\right) s^*=\mathbf{0}_{Nd}, \label{lem6 1}\\
& \mathcal{L} u^*=\mathbf{0}_{Nd},\quad \mathcal{L} s^*=\mathbf{0}_{Nd} ,\label{lem6 2}
\end{align}
where $\mathcal{L}=L\otimes I_d$ and $L$ is Laplacian matrix of graph $\mathcal{G}$. Due to  the properties of the Laplacian matrix, it is easy to derive that $u_i^*=u_j^*$ and $s_i^*=s_j^*$ for all $i\neq j$. Because of formulas \eqref{lem6 1}-\eqref{lem6 2}, it leads to
\begin{align}
    u_i^* & =\frac{1}{N} \sum_{i=1}^N \phi_i\left(x_{i}^*\right)=u(x^*), \label{lem6 3}\\
s_i^* & =\frac{1}{N} \sum_{i=1}^N \nabla_2 f_i\left(x_{i}^*, u(x^*)\right) .\label{lem6 4}
\end{align}
{By substituting \eqref{lem6 3} and \eqref{lem6 4} into \eqref{lem6 1}, we can obtain $\nabla F(x^*)=0$. Because $F$ is $\mu$-strongly convex, $x^*$ is the unique optimal solution to problem (1).}
\end{proof}
\vspace{-1em}
\subsection{Auxiliary Results}
In order to analyze the convergence and convergence rate of this algorithm, we use the method of compressed state vector and collect the following four quantities:
\begin{enumerate}[1)]
    \item $||x_{k+1}-x^*||$, the state error in the network;
    \item $||x_{k+1}-x_{k}||$, the state difference;
    \item $||u_{k+1}-\mathcal{K}u_{k+1}||$, the aggregative variable tracking error;
    \item $||s_{k+1}-\mathcal{K}s_{k+1}||$,  the gradient sum  tracking error.
\end{enumerate}
In the next Lemmas 7–10, we derive the relationships among the four quantities mentioned above. Firstly, we derive the bound on  $||x_{k+1}-x^*||$, the state error in the network.

\begin{lemma}
    Under Assumptions $1$-$4$, the  following inequality holds, $\forall k\geq  0$:
\begin{equation}
    \begin{aligned}
        & \left\|x_{k+1}-x^*\right\| \\
        & \leq(1-\mu \alpha)\left\|x_k-x^*\right\|+\alpha L_1\left\|u_k-\mathcal{K} u_k\right\|\\
&\quad+\alpha L_3\left\|s_k-\mathcal{K} s_k\right\|+\beta\left\|x_k-x_{k-1}\right\|.
    \end{aligned}
\end{equation}
\end{lemma}
\begin{proof}
For $\left\|x_{k+1}-x^*\right\|$, by invoking \eqref{compact hx}, it leads to
        \begin{align}
            & \left\|x_{k+1}-x^*\right\|\notag \\
& =\left\|x_k-x^*-\alpha\left[\nabla_1 f\left(x_k, u_k\right)+\nabla \phi\left(x_k\right) s_k\right]+\beta(x_k-x_{k-1})\right\| \notag \\
& \leq \left\|x_k-x^*-\alpha\left[\nabla_1 f\left(x_k, u_k\right)+\nabla \phi\left(x_k\right) s_k\right]\right\|+\beta\left\|x_k-x_{k-1}\right\| \notag \\
& \leq \Big\Vert x_k-x^*-\alpha\left[\nabla_1 f\left(x_k, \mathbf{1}_N \otimes \bar{u}_k\right)\right. \notag \\
&\quad \left.+\nabla \phi\left(x_k\right) \mathbf{1}_N \otimes \frac{1}{N} \sum_{i=1}^N \nabla_2 f_i\left(x_{i, k},  \bar{u}_k\right)\right]+\alpha \nabla F\left(x^*\right) \Big\Vert \notag \\
&\quad +\alpha \Big\Vert \nabla_1 f\left(x_k, u_k\right)+\nabla \phi\left(x_k\right) \mathbf{1}_N \otimes \bar{s}_k-\nabla_1 f\left(x_k, \mathbf{1}_N \otimes \bar{u}_k\right) \notag \\
&\quad -\nabla \phi\left(x_k\right) \mathbf{1}_N \otimes \frac{1}{N} \sum_{i=1}^N \nabla_2 f_i\left(x_{i, k}, \bar{u}_k\right) \Big\Vert \notag \\
&\quad +\alpha\left\|\nabla \phi\left(x_k\right) s_k-\nabla \phi\left(x_k\right) \mathbf{1}_N \otimes \bar{s}_k\right\|+\beta\left\|x_k-x_{k-1}\right\|. 
\label{7 big}
        \end{align}
From Lemma \ref{lem1}, we can bound the first term of the right term of \eqref{7 big} as follows:
    \begin{align}
        &\Big\Vert x_k-\alpha[\nabla \phi\left(x_k\right)
         \mathbf{1}_N \otimes \frac{1}{N} \sum_{i=1}^N \nabla_2 f_i\left(x_{i, k},  \bar{u}_k\right)\notag \\
         &\quad+\nabla_1 f\left(x_k, \mathbf{1}_N \otimes \bar{u}_k\right)]-\left[x^*-\alpha \nabla F\left(x^*\right)\right]\Big\Vert \notag\\
         &  \leq(1-\mu \alpha)\left\|x_k-x^*\right\|.
         \label{7 big 1}
    \end{align}
For the second term, since { $\nabla_1 f(x, {y})+\nabla \phi(x) 1_{{N}} \otimes \frac{1}{{N}} \sum_{i=1}^{{N}} \nabla_2 f_i\left(x_i, y_i\right)$} is  $L_1$-Lipschitz continuous  according to Assumption 2 and $\mathbf{1}_N \otimes \bar{u}_k=\mathcal{K} u_k$,  we can obtain
\begin{equation}
    \begin{aligned}
        & \alpha \Big\Vert \nabla_1 f\left(x_k, u_k\right)+\nabla \phi\left(x_k\right) \mathbf{1}_N \otimes \bar{s}_k-\nabla_1 f\left(x_k, \mathbf{1}_N \otimes \bar{u}_k\right) \\
& \quad-\nabla \phi\left(x_k\right) \mathbf{1}_N \otimes \frac{1}{N} \sum_{i=1}^N \nabla_2 f_i\left(x_{i, k},  \bar{u}_k\right) \Big\Vert \\
& \leq \alpha L_1\left\|u_k-\mathcal{K} u_k\right\|.
    \end{aligned}
\end{equation}
For the third term, using Assumption $4$ and $\mathbf{1}_N \otimes \bar{s}_k=\mathcal{K} s_k$ can get the following inequality:
\begin{equation}
    \begin{aligned}
        \alpha\left\|\nabla \phi\left(x_k\right) s_k-\nabla \phi\left(x_k\right) \mathbf{1}_N \otimes \bar{s}_k\right\| \leq \alpha L_3\left\|s_k-\mathcal{K} s_k\right\|.
    \end{aligned}
    \label{7 big 2}
\end{equation}
Then by substituting \eqref{7 big 1}-\eqref{7 big 2} into \eqref{7 big}, we complete the proof.
\end{proof}

Secondly, we derive the bound of $\left\|x_{k+1}-x_k\right\|$.

\begin{lemma}
    Under Assumptions $1$-$4$, the  following inequality holds, $\forall k\geq  0$:
\begin{equation}
    \begin{aligned}
\left\|x_{k+1}-x_k\right\| \leq & \alpha L_1\left(1+L_3\right)\left\|x_k-x^*\right\|+\alpha L_1\left\|u_k-\mathcal{K} u_k\right\| \\
& +\alpha L_3\left\|s_k-\mathcal{K} s_k\right\|+\beta\left\|x_k-x_{k-1}\right\|.
\end{aligned}
\label{lem 8}
\end{equation}
\end{lemma}
\begin{proof}
Note that $\nabla F(x^*)=0$ and then we have
    \begin{equation}
        \begin{aligned}
            & \left\|x_{k+1}-x_k\right\| \\
& =\left\|\beta(x_k-x_{k-1})-\alpha(\nabla_1 f\left(x_k, u_k\right)+\nabla \phi\left(x_k\right) s_k)\right\| \\
& \leq \alpha \Bigg\Vert \nabla_1 f\left(x_k, u_k\right)+\nabla \phi\left(x_k\right) \mathcal{K} s_k-\nabla_1 f\left(x^*, \mathbf{1}_N \otimes \bar{u}^*\right) \\
& \quad\quad-\nabla \phi\left(x^*\right)\left[\mathbf{1}_N \otimes \frac{1}{N} \sum_{i=1}^N \nabla_2 f_i\left(x^*_i,  \bar{u}^*\right)\right] \Bigg\Vert \\
& \quad+\alpha\left\|\nabla \phi\left(x_k\right)\left(s_k-\mathcal{K} s_k\right)\right\|+\beta\left\|x_k-x_{k-1}\right\|, \\
        \end{aligned}
        \label{8 big}
    \end{equation}
where $\bar{u}^*=\frac{\mathbf{1}_{N}^{\top}}{N}\otimes I_{d}u^* \in \mathbb{R}^d$. By utilizing Assumption $2$ and triangle inequality of norm, we can obtain the following formula:
{\begin{equation}
    \begin{aligned}
       &   \Big\Vert \nabla_1 f\left(x_k, u_k\right)+\nabla \phi\left(x_k\right) \mathcal{K} s_k-\nabla_1 f\left(x^*, \mathbf{1}_N \otimes \bar{u}^*\right) \\
& \quad-\nabla \phi\left(x^*\right)\left[\mathbf{1}_N \otimes \frac{1}{N} \sum_{i=1}^N \nabla_2 f_i\left(x^*_i,  \bar{u}^*\right)\right] \Big\Vert \\ 
& \leq L_1\left(\left\|x_k-x^*\right\|+\left\|u_k-\mathbf{1}_N \otimes \bar{u}^*\right\|\right)\\
& \leq L_1\left(\left\|x_k-x^*\right\|+\left\|u_k-\mathcal{K}u_k\right\|+\left\|\mathcal{K}u_k-\mathbf{1}_N \otimes \bar{u}^*\right\|\right).
    \end{aligned}
    \label{8 big1}
\end{equation}}
For $\left\|\mathcal{K}u_k-\mathbf{1}_N \otimes \bar{u}^*\right\|$, we can derive
    \begin{align}
\left\|\mathcal{K} u_k-\mathbf{1}_N \otimes \bar{u}^*\right\|^2 & =\left\|\mathbf{1}_N \otimes\left(\bar{u}_k-\bar{u}^*\right)\right\|^2 \notag\\
& =N\left\|\frac{1}{N} \sum_{i=1}^N\left(\phi_i\left(x_{i, k}\right)-\phi_i\left(x_i^*\right)\right)\right\|^2 \notag\\
& \leq \frac{1}{N}\left(\sum_{i=1}^N\left\|\phi_i\left(x_{i, k}\right)-\phi_i\left(x_i^*\right)\right\|\right)^2 \notag\\
& \leq \frac{1}{N}\left(\sum_{i=1}^N L_3\left\|x_{i, k}-x_i^*\right\|\right)^2 \notag\\
& \leq L_3^2 \sum_{i=1}^N\left\|x_{i, k}-x_i^*\right\|^2 \notag\\
& =L_3^2\left\|x_k-x^*\right\|^2,
\end{align}
where using the property that $\phi_i$ is $L_3$-Lipschitz continuous can obtain the second inequality, and applying the fact that $\left(\sum_{i=1}^N a_i\right)^2 \leq$ $N \sum_{i=1}^N a_i^2$ for any nonnegative scalars $a_i$s  can easily get the last inequality. Thus, we have
\begin{equation}
    \left\|\mathcal{K} u_k-\mathbf{1}_N \otimes \bar{u}^*\right\| \leq L_3\left\|x_k-x^*\right\|.
    \label{8 big2}
\end{equation}
Then by using Assumption $4$,  we can obtain
\begin{equation}
    \left\|\nabla \phi\left(x_k\right)\left(s_k-\mathcal{K} s_k\right)\right\| \leq L_3 \left\|s_k-\mathcal{K} s_k\right\|.
    \label{8 big3}
\end{equation}
Finally, by inserting \eqref{8 big1}, \eqref{8 big2} and \eqref{8 big3} into \eqref{8 big} we can obtain the result \eqref{lem 8}.
\end{proof}

The next step is to bound  the aggregative variable tracking error $||u_{k+1}-\mathcal{K}u_{k+1}||$.

\begin{lemma}
    
Under Assumptions $1$-$4$, the  following inequality holds, $\forall k\geq  0$:
\begin{equation}
\begin{aligned}
&\left\|u_{k+1}-\mathcal{K} u_{k+1}\right\| \\
& \leq\left(\rho+\alpha L_1 L_3\right)\left\|u_k-\mathcal{K} u_k\right\|+\alpha L_1 L_3\left(1+L_3\right)\left\|x_k-x^*\right\| \\
&\quad+\alpha L_3^2\left\|s_k-\mathcal{K} s_k\right\|+\beta L_3\left\|x_k-x_{k-1}\right\|.
\end{aligned}
\label{lem 9}
\end{equation}
\end{lemma}

\begin{proof}
For $\left\|u_{k+1}-\mathcal{K} u_{k+1}\right\|$, by invoking \eqref{compact hu}, it leads to
    \begin{equation}
\begin{aligned}
& \left\|u_{k+1}-\mathcal{K} u_{k+1}\right\| \\
& =\left\|\mathcal{A} u_k+\phi\left(x_{k+1}\right)-\phi\left(x_k\right)-\mathcal{K} \mathcal{A} u_k-\mathcal{K}\left[\phi\left(x_{k+1}\right)-\phi\left(x_k\right)\right]\right\| \\
& \leq \rho\left\|u_k-\mathcal{K} u_k\right\|+\|I_{Nd}-\mathcal{K}\|\left\|\phi\left(x_{k+1}\right)-\phi\left(x_k\right)\right\| \\
& \leq \rho\left\|u_k-\mathcal{K} u_k\right\|+L_3\|I_{Nd}-\mathcal{K}\|\left\|x_{k+1}-x_k\right\|,
\end{aligned}
\label{9 big}
\end{equation}
where Lemma $2$ has been utilized to obtain the first inequality, and by using  Assumption $4$ we can obtain  the last inequality.  Notice that $\|I_{Nd}-\mathcal{K}\|=1$, then by substituting \eqref{lem 8} into \eqref{9 big} we can complete the proof.
\end{proof}

Lastly, we derive the bound  $||s_{k+1}-\mathcal{K}s_{k+1}||$,  the tracking
 error of gradient sum  $\frac{1}{N}\sum_{i=1}^N\nabla_2 f_i\left(x_{i}, u(x_{k+1})\right)$.

\begin{lemma}
    
 Under Assumptions $1$-$4$, the  following inequality holds, $\forall k\geq  0$:
\begin{equation}
    \begin{aligned}
        & \left\|s_{k+1}-\mathcal{K} s_{k+1}\right\| \\
& \leq \left(\rho+\alpha L_2 L_3\left(1+L_3\right)\right)\left\|s_k-\mathcal{K} s_k\right\| 
  \\
& \quad+\left(\alpha L_1 L_2\left(1+L_3\right)+2L_2\right)\left\|u_k-\mathcal{K} u_k\right\|  \\
&\quad+ \alpha L_1 L_2\left(1+L_3\right)^2\left\|x_k-x^*\right\|+\beta L_2(1+L_3)||x_k-x_{k-1}||.
    \end{aligned}
\end{equation}
\end{lemma}

\begin{proof}For $\left\|s_{k+1}-\mathcal{K} s_{k+1}\right\|$, by invoking \eqref{compact hs}, it leads to
    \begin{equation}
\begin{aligned}
& \left\|s_{k+1}-\mathcal{K} s_{k+1}\right\| \\
& =||\mathcal{A} s_k+\nabla_2 f\left(x_{k+1}, u_{k+1}\right)-\nabla_2 f\left(x_k, u_k\right)\\
&\quad\quad-\mathcal{K}\mathcal{A} s_k-\mathcal{K}[\nabla_2 f\left(x_{k+1}, u_{k+1}\right)-\nabla_2 f\left(x_k, u_k\right)]|| \\
& \leq \left\|\mathcal{A}s_k-\mathcal{K} s_k\right\|\\
& \quad+\|I_{Nd}-\mathcal{K}\|\left\|\nabla_2 f\left(x_{k+1}, u_{k+1}\right)-\nabla_2 f\left(x_k, u_k\right)\right\|\\
& \leq \rho\left\|s_k-\mathcal{K} s_k\right\|+\left\|\nabla_2 f\left(x_{k+1}, u_{k+1}\right)-\nabla_2 f\left(x_k, u_k\right)\right\| \\
& \leq \rho\left\|s_k-\mathcal{K} s_k\right\|+  L_2\left(\left\|x_{k+1}-x_k\right\|+\left\|u_{k+1}-u_k\right\|\right),\\
\end{aligned}
\label{10 big}
\end{equation}
where  Lemma $2$ has been utilized to obtain the first inequality and Assumption $3$ has been leveraged in the last inequality.  Notice that
\begin{align}
    & \left\|u_{k+1}-u_k\right\|\notag\\
    & =\left\|\mathcal{A}u_k-u_k+\phi\left(x_{k+1}\right)-\phi\left(x_k\right)\right\|\notag\\
    & =\left\|\mathcal{A}u_k-\mathcal{A}\mathcal{K}u_k+\mathcal{K}u_k-u_k+\phi\left(x_{k+1}\right)-\phi\left(x_k\right)\right\|\notag\\
    & \leq \|(\mathcal{A}-\left.I_N \otimes I_d\right)\left(u_k-\mathcal{K} u_k\right)\|+\|\phi\left(x_{k+1}\right)-\phi\left(x_k\right)\|\notag\\
    & \leq ||A-I_N|| \left\|u_k-\mathcal{K} u_k\right\|+ L_3\left\|x_{k+1}-x_k\right\| \notag\\
    & \leq 2 \left\|u_k-\mathcal{K} u_k\right\|+ L_3\left\|x_{k+1}-x_k\right\|. 
    \label{10 big1}
    \end{align}
Then by substituting  \eqref{10 big1} and \eqref{lem 8} into \eqref{10 big}, we can finish the proof.
\end{proof}
\vspace{-1em}
\subsection{Main Result}
We now present the main result of this section. Based on Lemmas $7$-$10$, we give the convergence and convergence rate of the DAGT-HB algorithm in the following theorem.

\begin{theorem}
    
Under Assumptions 1-4, if 
\begin{equation}
    (\alpha,\beta) \in \bigcap\limits _{i=1}^6\mathcal{S}_i,
\end{equation}
where $\mathcal{S}_i,i=1,\cdots,6$  are defined in the following proof, then $x_k=$ col$\left(x_{1, k}, \ldots, x_{N, k}\right)$ generated by DAGT-HB  {converges to the
optimizer of problem (1) at the linear convergence rate.}
\end{theorem}

\begin{proof}
 Denote
\begin{equation}
    V_k=\text{col}\left(\left\|x_k-x^*\right\|,\left\|x_k-x_{k-1}\right\|,\left\|u_k-\mathcal{K} u_k\right\|,\left\|s_k-\mathcal{K} s_k\right\|\right).
\end{equation}
From Lemmas $7$-$10$, it can be concluded that
\begin{equation}
    V_{k+1} \leq P V_k,
\end{equation}
where 
\begin{equation}
     \begin{aligned}
     P=\left[\begin{array}{cccc}
1-\mu\alpha & \beta & \alpha L_1 & \alpha L_3\\
\alpha L_1(1+L_3) & \beta &\alpha L_1 & \alpha L_3\\
\alpha L_1L_3(1+L_3) & \beta L_3 & \rho+\alpha L_1L_3 & \alpha L_3^2\\
\alpha L_1L_3(1+L_3)^2 & \beta L_2(1+L_3) & p_{43} & p_{44}
\end{array}\right],\\
     \end{aligned}
 \end{equation}
 $p_{43}=\alpha L_1L_3(1+L_3)+2L_2$ and $p_{44}=\rho+\alpha L_2L_3(1+L_3).$
Firstly, based on Lemma $3$, we seek a small range of $\alpha$ and $\beta$
 to satisfy $\rho(P)<1$. We define a positive vector $z=[z_1,z_2,z_3,z_4]^{\top} \in \mathbb{R}^4$ such that
 \begin{equation}
     Pz<z,
 \end{equation}
 which is equal to
 \begin{equation}
     \begin{aligned}
        & 0<\beta<\frac{\alpha(\mu z_1-L_1z_3-L_3z_4)}{z_2}=M_1,\\
       &  0<\beta<\frac{z_2-\alpha L_1(1+L_3)z_1-\alpha L_1z_3-\alpha L_3 z_4}{z_2}=M_2,\\
       &  0<\beta<\\
       &\frac{(1-\rho-\alpha L_1L_3)z_3-\alpha L_1L_3(1+L_3)z_1-\alpha L_3^2z_4}{L_3z_2}=M_3,\\
       & 0<\beta<\frac{[1-\rho-\alpha L_2L_3(1+L_3)]z_4-\alpha L_1L_2(1+L_3)^2z_1}{L_2(1+L_3)z_2},\\
       &- \frac{[\alpha L_1L_2(1+L_3)+2L_2]z_3}{L_2(1+L_3)z_2}=M_4.
     \end{aligned}
 \end{equation}
 {From the above inequalities, we derive
     \begin{align}
         & 0< \alpha <\frac{z_2}{L_1(1+L_3)z_1+L_1z_3+L_3 z_4}=J_1,\notag\\
         & 0< \alpha <\frac{1-\rho}{L_1L_3}=J_2,\notag\\
         & 0< \alpha <\frac{(1-\rho)z_3}{L_1L_3(1+L_3)z_1+L_1L_3z_3+L_3^2z_4}=J_3,\notag\\
         & 0< \alpha < \frac{1-\rho}{ L_2L_3(1+L_3)}=J_4,\notag\\
         & 0< \alpha < 
         \frac{(1-\rho)z_4-2L_2z_3}{L_2(1+L_3)[L_1(1+L_3)z_1+L_1z_3+L_3z_4]}=J_5,\notag\\
         & z_1>\frac{L_1z_3+L_3z_4}{\mu},\quad z_4>\frac{2L_2z_3}{1-\rho}.
     \end{align}
 That is to say, we can select arbitrary $z_2$ and $z_3$, when
 \begin{equation}
     \begin{aligned}
       &  0<\alpha<\bar{\alpha}=\min \{ J_1, J_2, J_3, J_4, J_5,\frac{1}{L_1} \},
     \end{aligned}
 \end{equation}}
 and 
 \begin{equation}
     \begin{aligned}
         & 0<\beta<\bar{\beta}=\min \{M_1, M_2, M_3, M_4\},
     \end{aligned}
 \end{equation}
 where $z_2>0, z_3>0, z_1>\frac{L_1z_3+L_3z_4}{\mu}$ and $z_4>\frac{2L_2z_3}{1-\rho}$, we have $\rho(P)<1$. 

 Next, we use the  Jury criterion to seek precise range of $\alpha$ and $\beta$ to meet $\rho(P)<1$. By computing, we can obtain the characteristic polynomial of $P$:
 \begin{equation}
     H(\lambda)=|\lambda I-P|=a_0+a_1\lambda+a_2\lambda^2+a_3\lambda^3+\lambda^4,
 \end{equation}
 where
 \begin{equation}
     \begin{aligned}
         a_0 &=\beta d_1 \rho(\rho+2 \alpha d_2),\\
a_1&=  \beta[-d_1 \rho+(d_1+\rho)\times
\left(-\rho+2 \alpha d_2\right)] 
 +(\mu \alpha-1) \rho\left(\rho+\alpha d_2\right)\\
&\quad-\alpha L_3 d_1 \times
[\alpha L_1\left(\rho+
 \alpha L_3 d_2\right)+\alpha L_3 d_3 ],\\
 a_2&=\beta(d_1+2\rho+\alpha d_2)+(1-\mu\alpha)(2\rho+\alpha d_2)+\rho(\rho+\alpha d_2)\\
 &\quad+\alpha L_3 d_1[L_1+L_3(1+L_3)]+L_3[\alpha L_1(\rho+\alpha d_2)+\alpha L_3 d_3],\\
 a_3&=-\beta+(\mu-2)\alpha-1-\alpha L_3[L_2(1+L_3)+L_1],
     \end{aligned}
 \end{equation}
 and $d_1=1-\mu\alpha-\alpha L_1(1+L_3)$, $d_2=L_3(1+L_3)(L_2-L_3)$ and $d_3=-2\alpha L_2+\rho(1+L_3)$. Then we can obtain:
 \begin{equation}
     \begin{aligned}
         &b_0=a_0^2-1, b_1=a_0a_1-a_3, b_2=a_0a_2-a_2, b_3=a_0a_3-a_1,\\
         &c_0=b_0^2-b_3^2,c_1=b_0b_1-b_2b_3,c_2=b_0b_2-b_1b_3.
     \end{aligned}
 \end{equation}
 Next we denote
 \begin{equation}
     \begin{aligned}
         &\mathcal{S}_1=\{(\alpha,\beta)|a_0+a_1+a_2+a_3+1>0\},\\
        &\mathcal{S}_2=\{(\alpha,\beta)|a_0-a_1+a_2-a_3+1>0\},\\
        &\mathcal{S}_3=\{(\alpha,\beta)||a_0|<1\},\quad\mathcal{S}_4=\{(\alpha,\beta)||b_0|>|b_3|\},\\
        &\mathcal{S}_5=\{(\alpha,\beta)||c_0|>|c_2|\},\quad\mathcal{S}_6=\{(\alpha,\beta)|\alpha>0,\beta>0\}.
     \end{aligned}
 \end{equation}
 Hence, according to Lemma $5$, when $(\alpha,\beta) \in \bigcap\limits _{i=1}^6\mathcal{S}_i$, the spectral radius of the matrix $P$ is less than 1. In view of the above analysis, we know $\bigcap\limits _{i=1}^6\mathcal{S}_i$ is non-empty.
 Finally, denote $\rho_1=\rho(P)$ and $0<\rho_1<1$. Based on Lemma 4, we can derive that for $\epsilon \in \left(0, \frac{1-\rho_1}{2}\right)$, there exist some matrix norm $\| |\cdot|\|$ such
that $ \| |P|\|\leq \rho_1+\epsilon$. Then we can obtain 
\begin{equation}
    \| |V_{k+1}|\| \leq \| |P|\| \| |V_{k}|\| \leq (\rho(P)+\epsilon) \| |V_{k}|\| <  \frac{1+\rho_1}{2} \| |V_{k}|\|.
\end{equation}
By recalling that all norms are equivalent on finite-dimensional vector spaces, there always exist $\tau_1>0$ and $\tau_2>0$ such that $\|\cdot\| \leq \tau_1\| |\cdot|\|$ and $\| |\cdot|\| \leq \tau_2 \|\cdot\|$. Then we can obtain
 \begin{equation}
     ||x_k-x^*||\leq \|V_k\|\leq \tau_1 \||V_k|\| \leq C_1\left(\frac{1+\rho_1}{2}\right)^k,
 \end{equation}
 where $C_1=\tau_1\||V_0|\|$. So DAGT-HB achieves {the linear convergence rate.} Then the proof is completed.
\end{proof}
{
\begin{remark}
   Here, we point out the set $ \bigcap\limits _{i=1}^6\mathcal{S}_i$ is non-empty. In fact,  the set $\{ (\alpha,\beta)|0<\alpha<\bar{\alpha},0<\beta<\bar{\beta}\} \subset \bigcap\limits _{i=1}^6\mathcal{S}_i$, because the range of $\alpha$ and $\beta$ provided by the set  $\bigcap\limits _{i=1}^6\mathcal{S}_i$ is precise by leveraging Jury
criterion. Hence, when the set $\{ (\alpha,\beta)|0<\alpha<\bar{\alpha},0<\beta<\bar{\beta}\}$ is non-empty, {the set $\bigcap\limits _{i=1}^6\mathcal{S}_i$ is also non-empty.} It is obvious that the set $\{ (\alpha,\beta)|0<\alpha<\bar{\alpha},0<\beta<\bar{\beta}\}$ is not empty. For instance, we can select $z_2=z_3=1, z_1=\frac{2L_1z_3+L_3z_4}{\mu}$ and $z_4=\frac{3L_2z_3}{1-\rho}$, then we can obtain ${J}_1,\ldots,{J}_5$, ${M}_1,\ldots,{M}_4$ and the set $\{ (\alpha,\beta)|0<\alpha<\bar{\alpha},0<\beta<\bar{\beta}\}$ which is non-empty. Further, it can be deduced that the set $ \bigcap\limits _{i=1}^6\mathcal{S}_i$ is non-empty.
\end{remark}
}
\begin{remark}
    
 In Theorem $1$, we have established a linear rate of DAGT-HB when the step-size ${\alpha}$, and the  momentum
parameter ${\beta}$ follow (37). But we acknowledge
that the theoretical bounds of ${\alpha}$ and ${\beta}$ in Theorem 1 are conservative duo to the employment of inequality relaxation
techniques. How to obtain theoretical boundaries and even optimal parameters will be considered in our future work. Furthermore, we present a practical guide for selecting the two steps $\alpha$ and $\beta$ as follows:\\
\textbf{Step 1:} $L_1, L_2, L_3$ and $\mu$ are determined by  the  property
of the cost function. Moreover, $\rho$ is determined by the communication topology.\\
         \textbf{Step 2:} By computing the characteristic polynomial of $P$ and based on Jury criterion,  the set $\mathcal{S}_1,\ldots,\mathcal{S}_6$ can be obtained. Furthermore,  the theoretical bounds of $\alpha$ and $\beta$ can be derived. \\
          \textbf{Step 3:} By considering the theoretical results from Step 2 and  the specific circumstances in the practical application, parameters $\alpha$ and $\beta$ can be fine-tuned appropriately to achieve the best performance.

\end{remark}

 { When $x_k$ is generated by DAGT-HB, based on Theorem 1, we can obtain $F(x_k)$ converges to $F(x^*)$ with a linear convergence rate, which is summarized in the following corollary.}
\begin{corollary}
    
 Under the same assumptions of
Theorem $1$, the following equality holds:
\begin{equation}
    |f(x_k)-f(x^*)|\leq \frac{L_1}{2}C_1^2 \left(\frac{1+\rho_1}{2}\right)^{2k}.
\end{equation}
\end{corollary}
\begin{proof}Because $F(x)$ is $L_1$-smooth and $\nabla F(x^*)=0$, we can obtain
\begin{equation}
    F(x_k)-F(x^*)\leq \frac{L_1}{2}||x_k-x^*||^2.
\end{equation}
By substituting (50) into (52), we complete the proof.
\end{proof}

As a summary of this section, the DAGT-HB method is formulated as the following Algorithm \ref{algorithm1}.

\begin{algorithm}
   { \caption{DAGT-HB}\label{algorithm1}
    \SetKwInOut{Input}{Input}\SetKwInOut{Output}{Output}
    \Input{initial point $x_{i,-1} $ and $x_{i,0} \in \mathbb{R}^{n_i}$, $u_{i,0}=\phi_i(x_{i,0})$ and $s_{i,0}=\nabla_2 f_i(x_{i,0},u_{i,0})$ for $i=1,2,\ldots,N$, let $  (\alpha,\beta) \in \bigcap\limits _{i=1}^6\mathcal{S}_i$, set $k=0$, select appropriate $\epsilon_1>0$ and $K_{max} \in \mathbb{N}^+$. }
    \Output{optimal $x^*$ and $F(x^*)$.}
    While
    {$k<K_{max}$ and $||\nabla F(x_k)||\geq \epsilon_1$ do}\\
    Iterate: Update for each $i \in \{1, 2, \ldots, N\}$:
    \begin{align*}
    x_{i, k+1} & =x_{i, k}-\alpha\left[\nabla_1 f_i\left(x_{i, k}, u_{i, k}\right)+\nabla \phi_i\left(x_{i, k}\right) s_{i, k}\right] \notag\\
& \quad +\beta (x_{i,k}-x_{i,k-1}), \\
u_{i, k+1} & =\sum_{j=1}^N a_{i j} u_{j, k}+\phi_i\left(x_{i, k+1}\right)-\phi_i\left(x_{i, k}\right), \\
s_{i, k+1} & =\sum_{j=1}^N a_{i j} s_{j, k}+\nabla_2 f_i\left(x_{i, k+1}, u_{i, k+1}\right)  \notag\\
&\quad  -\nabla_2 f_i\left(x_{i, k}, u_{i, k}\right). 
\end{align*}

   Update: $k=k+1$.\\}
\end{algorithm}

\vspace{-1em}
\section{DAGT-NES}
In addition to heavy ball, {the Nesterov’s algorithm \cite{nesterov2003introductory}} is also a well-known accelerated  method that can be combined with DAGT to solve problem (1). To this end, we propose the following DAGT-NES algorithm:
\begin{align}
  &  x_{i, k+1}  =y_{i, k}-\alpha\left[\nabla_1 f_i\left(y_{i, k}, u_{i, k}\right)+\nabla \phi_i\left(y_{i, k}\right) s_{i, k}\right], \\
 &y_{i,k+1} = x_{i,k+1} + \gamma(x_{i,k+1}-x_{i,k}),\\
&u_{i, k+1}  =\sum_{j=1}^N a_{i j} u_{j, k}+\phi_i\left(y_{i, k+1}\right)-\phi_i\left(y_{i, k}\right), \\
&s_{i, k+1} =\sum_{j=1}^N a_{i j} s_{j, k}+\nabla_2 f_i\left(y_{i, k+1}, u_{i, k+1}\right) -\nabla_2 f_i\left(y_{i, k}, u_{i, k}\right). 
\end{align}
Likewise, in DAGT-NES, $u_{i,k} \in \mathbb{R}^d$ is introduced for agent $i$ to track the average $u(y)$ and $s_{i,k} \in \mathbb{R}^d$ tracks the gradient sum $\frac{1}{N}\sum_{i=1}^N\nabla_2 f_i\left(y_{i}, u(y)\right)$.
 The DAGT-NES algorithm can be rewritten as the following compact form:
\begin{align}
    x_{k+1} & =y_k-\alpha\left[\nabla_1 f\left(y_k, u_k\right)+\nabla \phi\left(y_k\right) s_k\right],\label{nesx} \\
    y_{k+1} & = x_{k+1}+\gamma (x_{k+1}-x_k),\label{nesy}\\
u_{k+1} & =\mathcal{A} u_k+\phi\left(y_{k+1}\right)-\phi\left(y_k\right), \label{nesu}\\
s_{k+1} & =\mathcal{A} s_k+\nabla_2 f\left(y_{k+1}, u_{k+1}\right)-\nabla_2 f\left(y_k, u_k\right).\label{ness}
\end{align}
with $\mathcal{A}=A \otimes I_d$ as defined in Lemma $2$ and $y_k=$ col$\left(y_{1, k}, \ldots, y_{N, k}\right) \in \mathbb{R}^n$. The  notations for $x_k \in \mathbb{R}^n$, $u_k$ and $s_k\in \mathbb{R}^{Nd}$ are the same as DAGT-HB.

Note that if  we initialize $u$ and $s$ as $u_{i,0}=\phi_i(y_{i,0})$ and $s_{i,0}=\nabla_2 f_i(y_{i,0},u_{i,0})$ for $i=1,2,\ldots,N$, where $y_{i,0}$ is arbitrary, analogous to DAGT-HB, we can obtain
\begin{align}
    \bar{u}_k & =\frac{1}{N} \sum_{i=1}^N u_{i, k}=\frac{1}{N} \sum_{i=1}^N \phi_i\left(y_{i, k}\right), \label{mean u}\\
\bar{s}_k & =\frac{1}{N} \sum_{i=1}^N s_{i, k}=\frac{1}{N} \sum_{i=1}^N \nabla_2 f_i\left(y_{i, k}, u_{i, k}\right) .\label{mean s}
\end{align}

Next, we establish the equivalence of the optimal solution to the problem (1) and the fixed point of the DAGT-NES algorithm. 

\begin{lemma}
    Under Assumption $1$ and Assumption $2$, the equilibrium point of \eqref{nesx}-\eqref{ness} is the optimal solution to problem (1).
\end{lemma}

\begin{proof}See Appendix 8.1.
\end{proof}
\vspace{-1em}
\subsection{Auxiliary Results}
Similar to DAGT-HB, we utilize  the method of compressed state vector to derive the convergence and convergence rate of DAGT-NES and still collect the following four quantities:
\begin{enumerate}[1)]
    \item $||x_{k+1}-x^*||$, the state error in the network;
    \item $||x_{k+1}-x_{k}||$, the state difference;
    \item $||u_{k+1}-\mathcal{K}u_{k+1}||$, the aggregative variable tracking error;
    \item $||s_{k+1}-\mathcal{K}s_{k+1}||$,  the gradient sum  tracking error.
\end{enumerate}

In the next Lemmas 12–15, we derive the relationship among the four quantities mentioned above. Firstly, we derive the bound on  $||x_{k+1}-x^*||$, the state error in the network.

\begin{lemma}
    
Under Assumptions $1$-$4$, the  following inequality holds, $\forall k\geq  0$:
\begin{equation}
    \begin{aligned}
        & \left\|x_{k+1}-x^*\right\| \\
        & \leq(1-\mu \alpha)\left\|x_k-x^*\right\|+\alpha L_1\left\|u_k-\mathcal{K} u_k\right\|\\
&\quad+\alpha L_3\left\|s_k-\mathcal{K} s_k\right\|+(1-\mu \alpha)\gamma\left\|x_k-x_{k-1}\right\|.\\
    \end{aligned}
\end{equation}
\end{lemma}

\begin{proof}See Appendix 8.2.
\end{proof}

Secondly, we derive the bound of $\left\|x_{k+1}-x_k\right\|$.

\begin{lemma}
    
Under Assumptions 1-4, the  following inequality holds, $\forall k\geq  0$:
\begin{equation}
    \begin{aligned}
        &\left\|x_{k+1}-x_k\right\|\\
        &\leq  \alpha L_1\left(1+L_3\right)\left\|x_k-x^*\right\|+\alpha L_1\left\|u_k-\mathcal{K} u_k\right\| \\
& \quad+\alpha L_3\left\|s_k-\mathcal{K} s_k\right\|+\gamma(1+\alpha L_1+\alpha L_1L_3)\left\|x_k-x_{k-1}\right\|.
    \end{aligned}
    \label{lem14}
\end{equation}
\end{lemma}

\begin{proof}See Appendix 8.3.
\end{proof}

The next step is to bound  the aggregative variable tracking error $||u_{k+1}-\mathcal{K}u_{k+1}||$.

\begin{lemma}
    
Under Assumptions $1$-$4$, the  following inequality holds, $\forall k\geq  0$:
\begin{equation}
    \begin{aligned}
       &\left\|u_{k+1}-\mathcal{K} u_{k+1}\right\| \\
& \leq\left[\rho+\alpha L_1 L_3(\gamma+1)\right]\left\|u_k-\mathcal{K} u_k\right\|+\alpha L_3^2(\gamma+1)\left\|s_k-\mathcal{K} s_k\right\|\\
&\quad +\alpha L_1 L_3\left(1+L_3\right)(\gamma+1)\left\|x_k-x^*\right\| \\
&\quad +\gamma L_3[(1+\gamma)(1+\alpha L_1 +\alpha L_1L_3)+1]\left\|x_k-x_{k-1}\right\|.\\ 
    \end{aligned}
\end{equation}
\end{lemma}

\begin{proof}See Appendix 8.4.
\end{proof}

Lastly, we derive the bound  $||s_{k+1}-\mathcal{K}s_{k+1}||$,  the tracking
 error of gradient sum $\frac{1}{N}\sum_{i=1}^N\nabla_2 f_i\left(y_{i}, u(y_{k+1})\right)$.

\begin{lemma}
    
Under Assumptions $1$-$4$, the  following inequality holds, $\forall k\geq  0$:
\begin{equation}
    \begin{aligned}
        & \left\|s_{k+1}-\mathcal{K} s_{k+1}\right\| \\
& < \left[\rho+\alpha L_2 L_3\left(1+L_3\right)(1+\gamma)\right]\left\|s_k-\mathcal{K} s_k\right\| \\
&\quad +\alpha L_1 L_2\left(1+L_3\right)^2(\gamma+1)\left\|x_k-x^*\right\| \\
& \quad+\left[\alpha L_1 L_2\left(1+L_3\right)(r+1)+2L_2\right]\left\|u_k-\mathcal{K} u_k\right\|  \\
&\quad +\gamma L_2(L_3+1)[(1+\gamma)(1+\alpha L_1 +\alpha L_1L_3)+1]||x_k-x_{k-1}||.
    \end{aligned}
\end{equation}
\end{lemma}

\begin{proof}See Appendix 8.5.
\end{proof}

\vspace{-1em}

\subsection{Main Result}
By summarizing 
{Lemmas 12-15}, we give the convergence and convergence rate  of the DAGT-NES algorithm in the following theorem.

\begin{theorem}
    
 Under Assumptions 1-4, if 
\begin{equation}
   (\alpha,\gamma) \in \bigcap\limits _{i=1}^6\mathcal{T}_i,
\end{equation}
where $\mathcal{T}_i$, $i=1,\cdots,6$  are defined in the following proof, then $x_k=$ col$\left(x_{1, k}, \ldots, x_{N, k}\right)$ generated by DAGT-NES {converges to the
optimizer of problem (1) at the linear convergence rate.}
\end{theorem}
\begin{proof}
    Denote
\begin{equation}
    V_k=\text{col}\left(\left\|x_k-x^*\right\|,\left\|x_k-x_{k-1}\right\|,\left\|u_k-\mathcal{K} u_k\right\|,\left\|s_k-\mathcal{K} s_k\right\|\right).
\end{equation}
From Lemmas {12-15}, it can be concluded that
\begin{equation}
    V_{k+1} \leq Q V_k,
\end{equation}
where 
\begin{equation}
     \begin{aligned}
     Q=\left[\begin{array}{cccc}
1-\mu\alpha & (1-\mu\alpha)\gamma & \alpha L_1 & \alpha L_3\\
\alpha L_1(1+L_3) & q_{22} &\alpha L_1 & \alpha L_3\\
\alpha L_1 L_3\left(1+L_3\right)(\gamma+1) & q_{32} & q_{33} & q_{34} \\
\alpha L_1L_2(1+L_3)^2 (\gamma+1)& q_{42} & q_{43} & q_{44}
\end{array}\right],\\
     \end{aligned}
 \end{equation}
and 
\begin{equation}
\left\{\begin{array}{l}
 q_{22}=\gamma(1+\alpha L_1+\alpha L_1L_3),\\
 q_{32}=\gamma L_3[(1+\gamma)(1+\alpha L_1 +\alpha L_1L_3)+1],\\
q_{33}=\rho+\alpha L_1 L_3(\gamma+1),\\
q_{34}=\alpha L_3^2(\gamma+1),\\
q_{42}=\gamma L_2(L_3+1)[(1+\gamma)(1+\alpha L_1 +\alpha L_1L_3)+1],\\
q_{43}=\alpha L_1 L_2\left(1+L_3\right)(\gamma+1)+2L_2,\\
q_{44}=\rho+\alpha L_2 L_3\left(1+L_3\right)(1+\gamma).
\end{array}\right.
\end{equation}
Firstly, based on Lemma $3$, we want to seek a small range of $\alpha$ and $\gamma$
 to satisfy $\rho(Q)<1$. We define a positive vector $t=[t_1,t_2,t_3,t_4]^{\top}$ such that
 \begin{equation}
     Qt<t .
 \end{equation}
However, since in the matrix $Q$, $\alpha$ and $\gamma$ have strong nonlinear relationship,  it is difficult to give the  range of $\alpha$ and $\gamma$ that makes the spectral radius of the matrix $Q$ less than 1. Thus,  to simplify the calculation, we let 
\begin{equation}
    \alpha\leq \frac{1}{L_1} \quad \text{and} \quad \gamma\leq \min \{\frac{1}{L_2},\frac{1}{L_3}\}
\end{equation}
 to eliminate some entries in matrix $Q$ that contain the nonlinear relationship between $\alpha$ and $\gamma$. Then we can obtain 
 \begin{equation}
    V_{k+1}<R V_k,
\end{equation}
where 
\begin{equation}
     \begin{aligned}
     R=\left[\begin{array}{cccc}
1-\mu\alpha & (1-\mu\alpha)\gamma & \alpha L_1 & \alpha L_3\\
\alpha L_1(1+L_3) & \gamma(2+L_3) &\alpha L_1 & \alpha L_3\\
\alpha L_1 L_3\left(2+L_3\right) & r_{32} & r_{33} & \alpha L_3(L_3+1)\\
r_{41} & r_{42} & r_{43} & r_{44}
\end{array}\right],\\
     \end{aligned}
 \end{equation}
and
\begin{equation}
\left\{\begin{array}{l}
 r_{32}=\gamma(L_3^2+4L_3+2),\\
r_{33}=\rho+\alpha L_1(L_3+1),\\
r_{41}=\alpha L_1(1+L_2)(1+L_3)^2,\\
r_{42}=\gamma(L_3+1)(L_2L_3+2L_2+L_3+1),\\
r_{43}=\alpha L_1 (L_2+1)\left(1+L_3\right)+2L_2,\\
r_{44}=\rho+\alpha L_2 \left(1+L_3\right)^2.
\end{array}\right.
\end{equation}
{Then we  solve the following equality:}
\begin{equation}
    Rt<t,
\end{equation}
which is equivalent to
\begin{equation}
     \begin{aligned}
        & 0<\gamma<\frac{\alpha(\mu t_1-L_1t_3-L_3t_4)}{(1-\mu\alpha)t_2}=\Gamma_1,\\
       &  0<\gamma<\frac{t_2-\alpha L_1(1+L_3)t_1-\alpha L_1t_3-\alpha L_3 t_4}{(2+L_3)t_2}=\Gamma_2,\\
       &  0<\gamma<\frac{[1-\rho-\alpha L_1(L_3+1)]t_3-\alpha L_1L_3(2+L_3)t_1}{(L_3^2+4L_3+2)t_2}\\
       & -\frac{\alpha L_3(L_3+1)t_4}{(L_3^2+4L_3+2)t_2}=\Gamma_3,\\
       & 0<\gamma<\frac{[1-\rho-\alpha L_2(1+L_3)^2]t_4-\alpha L_1(L_2+1)(1+L_3)^2t_1}{L_2(1+L_3)t_2}\\
       &- \frac{[\alpha L_1(L_2+1)(1+L_3)+2L_2]t_3}{L_2(1+L_3)t_2}=\Gamma_4.
     \end{aligned}
 \end{equation}
 From the above inequalities, we derive
     \begin{align}
         & 0< \alpha <\frac{t_2}{L_1(1+L_3)t_1+L_1t_3+L_3 t_4}=\Theta_1,\notag\\
         & 0< \alpha <\frac{1-\rho}{L_1(L_3+1)}=\Theta_2,\notag\\
         & 0< \alpha <\frac{(1-\rho)t_3}{L_1L_3(2+L_3)t_1+L_1(L_3+1)t_3+(L_3^2+L_3)t_4}=\Theta_3,\notag\\
         & 0< \alpha < \frac{1-\rho}{ L_2L_3(1+L_3)}=\Theta_4,\notag\\
         & 0< \alpha < \notag\\
         &\frac{(1-\rho)t_4-2L_2t_3}{L_1(L_2+1)(1+L_3)[(1+L_3)t_1+t_3]+L_2(1+L_3)^2t_4}=\Theta_5,\notag\\
         & t_1>\frac{L_1t_3+L_3t_4}{\mu}, \quad t_4>\frac{2L_2t_3}{1-\rho}.
     \end{align}
 That is to say, we can select arbitrary $t_2$ and $t_3$, when
 \begin{equation}
     \begin{aligned}
       &  0<\alpha<\bar{\alpha}=\min \{\Theta_1,\Theta_2,\Theta_3,\Theta_4,\Theta_5,\frac{1}{L_1}  \},
     \end{aligned}
 \end{equation}
 and 
 \begin{equation}
     \begin{aligned}
         & 0<\gamma<\bar{\gamma}=\min \{\Gamma_1,\Gamma_2,\Gamma_3,\Gamma_4,\frac{1}{L_2},\frac{1}{L_3}\},
     \end{aligned}
 \end{equation}
 where $t_2>0$, $t_3>0$, $t_1>\frac{L_1t_3+L_3t_4}{\mu}$ and $t_4>\frac{2L_2t_3}{1-\rho}$, we have $\rho(R)<1$. 

Next, we use Jury criterion to seek precise range of $\alpha$ and $\beta$ to meet $\rho(Q)<1$. By computing, we can obtain the characteristic polynomial of $Q$:
 \begin{equation}
     G(\lambda)=|\lambda I-Q|=a_0+a_1\lambda+a_2\lambda^2+a_3\lambda^3+\lambda^4,
 \end{equation}
 where
 \begin{equation}
     \begin{aligned}
         a_0&=(e_1-1)[\alpha L_3(\rho\alpha L_1-e_2)+\rho^2e_3]+\rho^2\alpha\gamma e_1L_1(1+L_3),\\
         a_1&=(e_1-1)[L_3e_2(1+r)-e_4\gamma+\rho(2e_3+\rho-\alpha L_3^2)]\\
         &\quad+\gamma L_3[e_2+\alpha L_3(e_1+\rho e_1-1-2\rho)]\\
         &\quad -\rho^2e_3-\alpha L_1(1+L_3)\rho (2\gamma e_1+\rho),\\
         a_2&=(1+\gamma)L_3e_2+(e_1-1)[(1+\gamma)e_4+2\rho+e_3]\\
         &\quad+\alpha L_1[(1+L_3)(2\rho+\gamma e_1)+L_3\gamma+L_3(1+\gamma)(e_1-1-\rho)]\\
         &\quad -\gamma e_4+\rho(2e_3+\rho),\\
         a_3&=e_1-e_2-1-2\rho-(1+\gamma)e_4-\alpha L_1[L_3(2+\gamma)+1] ,
         \end{aligned}
 \end{equation}
 and $e_1=\alpha[\mu+L_1(1+L_3)]$, $e_2=\alpha L_2[\rho(1+L_3)-2L_3]$, $e_3=\gamma(1+\alpha L_1+\alpha L_1L_3)$ and $e_4=\alpha L_2L_3(1+L_3)$. Then we can obtain:
 \begin{equation}
     \begin{aligned}
         &b_0=a_0^2-1, b_1=a_0a_1-a_3, b_2=a_0a_2-a_2, b_3=a_0a_3-a_1,\\
         &c_0=b_0^2-b_3^2,c_1=b_0b_1-b_2b_3,c_2=b_0b_2-b_1b_3.
     \end{aligned}
 \end{equation}
 Next we denote
 \begin{equation}
     \begin{aligned}
         &\mathcal{T}_1=\{(\alpha,\gamma)|a_0+a_1+a_2+a_3+1>0\},\\
        &\mathcal{T}_2=\{(\alpha,\gamma)|a_0-a_1+a_2-a_3+1>0\},\\
        &\mathcal{T}_3=\{(\alpha,\gamma)||a_0|<1\},\quad\mathcal{T}_4=\{(\alpha,\gamma)||b_0|>|b_3|\},\\
        &\mathcal{T}_5=\{(\alpha,\gamma)||c_0|>|c_2|\},\quad\mathcal{T}_6=\{(\alpha,\gamma)|\alpha>0,\gamma>0\}.
     \end{aligned}
 \end{equation}
 Thus, according to Lemma $5$, when $(\alpha,\gamma) \in \bigcap\limits _{i=1}^6\mathcal{T}_i$, the spectral radius of the matrix $Q$ is less than 1. In view of the above analysis, we know $\bigcap\limits _{i=1}^6\mathcal{T}_i$ is non-empty. Finally, denote $\rho_2=\rho(Q)$ and $0<\rho_2<1$. Based on Lemma 4, we can derive that for $\epsilon \in \left(0, \frac{1-\rho_2}{2}\right)$, there exist some matrix norm $\| |\cdot|\|$ such
that $ \| |Q|\|\leq \rho_2+\epsilon$. Then we can obtain 
\begin{equation}
    \| |V_{k+1}|\| \leq \| |Q|\| \| |V_{k}|\| \leq (\rho_2+\epsilon) \| |V_{k}|\| <  \frac{1+\rho_2}{2} \| |V_{k}|\|.
\end{equation}
By recalling that all norms are equivalent on finite-dimensional vector spaces, there always exist $\tau_3>0$ and $\tau_4>0$ such that $\|\cdot\| \leq \tau_3\| |\cdot|\|$ and $\| |\cdot|\| \leq \tau_4 \|\cdot\|$. Then we can obtain
\begin{equation}
     ||x_k-x^*||\leq \|V_k\|\leq \tau_3 \||V_k|\| \leq C_2\left(\frac{1+\rho_2}{2}\right)^k,
 \end{equation}
 where $C_2=\tau_3\||V_0|\|$. Thus, DAGT-NES achieves {the linear convergence rate.} Then the proof is completed.
\end{proof}

{
\begin{remark}
   Here, we point out the set $ \bigcap\limits _{i=1}^6\mathcal{T}_i$ is non-empty. In fact,  the set $\{ (\alpha,\gamma)|0<\alpha<\bar{\alpha},0<\gamma<\bar{\gamma}\} \subset \bigcap\limits _{i=1}^6\mathcal{T}_i$, because the range of $\alpha$ and $\gamma$ provided by the set  $\bigcap\limits _{i=1}^6\mathcal{T}_i$ is precise by leveraging Jury
criterion. Hence, when the set $\{ (\alpha,\gamma)|0<\alpha<\bar{\alpha},0<\gamma<\bar{\gamma}\}$ is non-empty, {the set $\bigcap\limits _{i=1}^6\mathcal{T}_i$ is also non-empty.} It is obvious that the set $\{ (\alpha,\gamma)|0<\alpha<\bar{\alpha},0<\gamma<\bar{\gamma}\}$ is non-empty. For instance, we can select $t_2=t_3=1, t_1=\frac{2L_1t_3+L_3t_4}{\mu}$ and $t_4=\frac{3L_2t_3}{1-\rho}$, then we can obtain $\Theta_1,\ldots,\Theta_5$, $\Gamma_1,\ldots,\Gamma_4$ and the set $\{ (\alpha,\gamma)|0<\alpha<\bar{\alpha},0<\gamma<\bar{\gamma}\}$ which is non-empty. Further, it can be deduced that the set $ \bigcap\limits _{i=1}^6\mathcal{T}_i$ is non-empty.
\end{remark}
}

\begin{remark}

 In Theorem $2$, we have established a linear rate of DAGT-NES when the  step-size ${\alpha}$, and the largest momentum
parameter ${\gamma}$ follow  (67). {Similar to DAGT-HB, the theoretical bounds of ${\alpha}$ and ${\gamma}$ in Theorem 2 are conservative duo to the employment of inequality relaxation
techniques. How to obtain theoretical boundaries and even optimal parameters will be considered in our future work. Furthermore, we present a practical guide for selecting the two steps $\alpha$ and $\gamma$ as follows:\\
\textbf{Step 1:} $L_1, L_2, L_3$ and $\mu$ are determined by  the  property
of the cost function. Moreover, $\rho$ is determined by the communication topology.\\
         \textbf{Step 2:} By computing the characteristic polynomial of $Q$ and based on Jury criterion,  the set $\mathcal{T}_1,\ldots,\mathcal{T}_6$ can be obtained. Furthermore,  the theoretical bounds of $\alpha$ and $\gamma$ can be derived. \\
          \textbf{Step 3:} By considering the theoretical results from Step 2 and  the specific circumstances in the practical application, parameters $\alpha$ and $\gamma$ can be fine-tuned appropriately to achieve the best performance.}
\end{remark}

\begin{remark}
Although our theoretical results are based on the assumption that communication graph is undirected and balanced, the DAGT-HB and DAGT-NES algorithms can be effectively extended to directed unbalanced graphs. This is achieved by introducing an auxiliary variable to estimate the Perron eigenvector of the weight matrix and establishing the consensus law \cite{8387465}. The relevant theoretical analysis will 
be considered in our future work. 
\end{remark}

 { When $x_k$ is generated by DAGT-NES, based on Theorem 1, we can obtain $F(x_k)$ converges to $F(x^*)$ with a linear convergence rate, which is summarized in the following corollary.}

{\begin{corollary}
    
Under the same assumptions of
Theorem 2, the following equality holds:
\begin{equation}
    |F(x_k)-F(x^*)|\leq \frac{L_1}{2}C_2^2 \left(\frac{1+\rho_2}{2}\right)^{2k}.
\end{equation}
\end{corollary}
}
\begin{proof}Because $F(x)$ is $L_1$-smooth and $\nabla F(x^*)=0$, we can obtain
\begin{equation}
    F(x_k)-F(x^*)\leq \frac{L_1}{2}||x_k-x^*||^2.
\end{equation}
By substituting (87) into (89), we complete the proof.
\end{proof}

As a summary of this section, the DAGT-NES method is formulated as the following Algorithm \ref{algorithm2}.

\begin{algorithm}
   { \caption{DAGT-NES}\label{algorithm2}
    \SetKwInOut{Input}{Input}\SetKwInOut{Output}{Output}
    \Input{initial point $x_{i,0} $ and $y_{i,0} \in \mathbb{R}^{n_i}$, $u_{i,0}=\phi_i(y_{i,0})$ and $s_{i,0}=\nabla_2 f_i(y_{i,0},u_{i,0})$ for $i=1,2,\ldots,N$, let $  (\alpha,\gamma) \in \bigcap\limits _{i=1}^6\mathcal{T}_i$, set $k=0$, select appropriate $\epsilon_2>0$ and $K_{max} \in \mathbb{N}^+$ .}
    \Output{optimal $x^*$ and $F(x^*)$.}
    While
    {$k<K_{max}$ and $||\nabla F(x_k)||\geq \epsilon_2$ do}\\
    Iterate: Update for each $i \in \{1, 2, \ldots, N\}$:
    \begin{align*}
    x_{i, k+1}  &=y_{i, k}-\alpha\left[\nabla_1 f_i\left(y_{i, k}, u_{i, k}\right)+\nabla \phi_i\left(y_{i, k}\right) s_{i, k}\right], \\
 y_{i,k+1}& = x_{i,k+1} + \gamma(x_{i,k+1}-x_{i,k}),\\
u_{i, k+1}  &=\sum_{j=1}^N a_{i j} u_{j, k}+\phi_i\left(y_{i, k+1}\right)-\phi_i\left(y_{i, k}\right), \\
s_{i, k+1} &=\sum_{j=1}^N a_{i j} s_{j, k}+\nabla_2 f_i\left(y_{i, k+1}, u_{i, k+1}\right)\\
&\quad-\nabla_2 f_i\left(y_{i, k}, u_{i, k}\right). 
\end{align*}

Update: $k=k+1$.\\}
    
\end{algorithm}
\vspace{-1 em}
\section{Discussion}
{
In Sections 3 and 4, we have demonstrated that the DAGT-HB and DAGT-NES algorithms can converge to the optimal solution at a global linear convergence rate under Assumptions 1-4. However, it is important to note that the linear convergence rate is established in \cite{Li2022Aggregative}. This leads to an intriguing question: what impact does the inclusion of the momentum term have on the algorithm’s convergence rate?

Here, we consider the following distributed quadratic optimization as an illustrative example:
\begin{equation}
\begin{aligned}
     &\min _{x \in \mathbb{R}^{N}} f(x), f(x)  =\sum_{i=1}^N f_i\left(x_i, u(x)\right), u(x)=\frac{\sum_{i=1}^N \phi_i(x_i)}{N},x_i\in \mathbb{R}\\
&f_i(x_i,u(x))=\frac{c_i}{2}x_i^{2}+\frac{u(x)}{N}, \phi_i(x_i)=h_ix_i+l_i, i=1,\ldots,N,
     \end{aligned}
     \label{example}
\end{equation}
where $c_i>0, h_i\geq 0$.  Then we can obtain $\mu=\min \{c_1,\ldots,c_N\}$ and $L_1=\max \{c_1,\ldots,c_N\}$. The quadratic optimization has been widely noticed in machine
 learning and also been considered in \cite{wang2024momentum}\cite{chen2024achieving}.

Let $\psi^d_k=$col$(x_k-x^*,u_k-K_Nu_k,s_k-K_Ns_k)$. Then by leveraging DAGT in \cite{Li2022Aggregative}, we can obtain that
 \begin{equation}
     \begin{aligned}
         \psi^d_{k+1}=P_D \psi^d_k+\theta_D
\end{aligned}
\label{pd}
 \end{equation}
 where
 $$
P_D=\left[\begin{array}{ccc}
I_N-\alpha C &  0 & 0 \\
 K_NH-H& A-K_N & 0\\
 0 & 0 & A-K_N
\end{array}\right] \in \mathbb{R}^{3N\times 3N},
 $$
 and 
 $$
\theta_D=\left[\begin{array}{c}
-\alpha H s_{k} \\
(H-{K}_NH)(x_{k+1}-x^*)\\
0
\end{array}\right] \in \mathbb{R}^{3N\times N}
 $$
where $C=\operatorname{diag}(c_1,\ldots,c_N) \mathbb{R}^{N\times N}$ and $H=\operatorname{diag}(h_1,\ldots,h_N) \mathbb{R}^{N\times N}$. By using the formulas \eqref{mean hu} and \eqref{mean hs}, we can obtain
\begin{equation}
    \begin{aligned}
        s_k&=s_k-{K}_Ns_k+{K}_Ns_k=s_k-{K}_Ns_k+\frac{1}{N}.
    \end{aligned}
    \label{sk}
\end{equation}
By bonding \eqref{pd} with \eqref{sk}, we can obtain
\begin{equation}
     \begin{aligned}
         P_{D1}\psi^d_{k+1}\leq P_{D2}\psi^d_{k},
\end{aligned}
\label{PD1}
 \end{equation}
where
 $$
P_{D1}=\left[\begin{array}{ccc}
I_N &  0 & 0 \\
 K_NH-H& I_N & 0\\
 0 & 0 & I_N
\end{array}\right] \in \mathbb{R}^{3N\times 3N},
 $$
 and 
  $$
P_{D2}=\left[\begin{array}{ccc}
I_N-\alpha C &  0 & -\alpha H \\
 K_NH-H& A-K_N & 0\\
 0 & 0 & A-K_N
\end{array}\right] \in \mathbb{R}^{3N\times 3N}.
 $$
Further, we can derive
\begin{equation}
    \psi^d_{k+1} \leq P_{D3} \psi^d_k,
\end{equation}
 where
 $$
P_{D3}=P_{D1}^{-1}P_{D2}=\left[\begin{array}{ccc}
I_N-\alpha C &  0 & -\alpha H \\
 \alpha(H-K_NH)C& A-K_N & 0\\
 0 & 0 & A-K_N
\end{array}\right].
 $$
 By employing Lemma 4, we can derive that there exists a matrix norm $\||\cdot|\|$ such that 
 \begin{align}
     \||\psi^d_{k+1}|\|\leq (\rho(P_{D3})+o(1)) \||\psi^d_{k}|\|.
 \end{align}
By  recalling that all norms are equivalent on finite-dimensional vector spaces, there always exist $\tau_5>0$ and $\tau_6>0$ such that $\|\cdot\| \leq \tau_5\| |\cdot|\|$ and $\| |\cdot|\| \leq \tau_6 \|\cdot\|$. Then we can obtain
 \begin{equation}
     ||x_k-x^*||\leq \|\psi^d_{k}\|\leq \tau_5 \||\psi^d_{k}|\| \leq \tau_5\left(\rho(P_{D3})+o(1)\right)^k\||\psi^d_{0}|\|.
 \end{equation}
 Hence, the convergence rate of DAGT is only related to the spectral radius of $P_{D3}$. By recalling the Lemma 2, we can obtain $\rho(P_{D3})=\max \{\rho(I_N-\alpha C),\rho\}$. Further, we can derive easily that $\rho(I_N-\alpha C)=\max \{|1-\alpha \mu|,|1-\alpha L_1|\}$. So when $\alpha=\frac{2}{\mu+L_1}$, $\rho(I_N-\alpha C)$ achieves the minimum $\frac{L_1-\mu}{L_1+\mu}$. Thus, for the quadratic optimization, the DAGT can achieve a linear convergence rate with $\mathcal{O}\left(\max\left\{\rho^k,(\frac{L_1-\mu}{L_1+\mu})^k\right\}\right)$. 

 Next, we establish the convergence rate of DAGT-HB for the quadratic optimization. Let  $\psi^h_k=$col$(x_k-x^*,x_{k-1}-x^*,u_k-K_Nu_k,s_k-K_Ns_k)$. Similar to the analysis of DAGT, we can obtain
 \begin{equation}
     \begin{aligned}
         \psi^h_{k+1}=P_H\psi^h_k+\theta_H,
\end{aligned}
\label{ph}
 \end{equation}
 where
 $$
P_H=\left[\begin{array}{cccc}
(1+\beta)I_N-\alpha C &-\beta I_N&  0 & 0 \\
I_N & 0 & 0 & 0\\
 K_NH-H&0 & A-K_N & 0\\
 0 & 0 &0 & A-K_N
\end{array}\right],
 $$
 and 
 $$
\theta_H=\left[\begin{array}{c}
-\alpha H s_{k} \\
0\\
(H-{K}_NH)(x_{k+1}-x^*)\\
0
\end{array}\right] .
 $$
 By employing \eqref{sk}, we have
\begin{align}
    P_{H1}\psi^h_{k+1} \leq P_{H2}\psi^h_{k},
    \label{PH1}
\end{align}
where
 $$
P_{H1}=\left[\begin{array}{cccc}
I_N & 0& 0 & 0 \\
0 &I_N & 0 &0\\
 K_NH-H&0& I_N & 0\\
 0& 0 & 0 & I_N
\end{array}\right],
 $$
 and 
 $$
P_{H2}=\left[\begin{array}{cccc}
(1+\beta)I_N-\alpha C &-\beta I_N&  0 & -\alpha H \\
I_N & 0 & 0 & 0\\
 K_NH-H&0 & A-K_N & 0\\
 0 & 0 &0 & A-K_N
\end{array}\right].
 $$
 Further, we can derive
 \begin{align}
     \psi^h_{k+1} \leq P_{H1}^{-1}P_{H2}\psi^h_{k}=P_{H3}\psi^h_{k},
 \end{align}
where
$$
P_{H3}=\left[\begin{array}{cccc}
(1+\beta)I_N-\alpha C &-\beta I_N&  0 & -\alpha I_N \\
I_N & 0 & 0 & 0\\
*&**& A-K_N &*** \\
 0 & 0 &0 & A-K_N
\end{array}\right].
 $$
 By employing Lemma 4, we can derive that there exists a matrix norm $\||\cdot|\|$ such that 
 \begin{align}
     \||\psi^h_{k+1}|\|\leq (\rho(P_{H3})+o(1)) \||\psi^h_{k}|\|.
 \end{align}
 By recalling that all norms are equivalent on finite-dimensional vector spaces, there always exist $\tau_7>0$ and $\tau_8>0$ such that $\|\cdot\| \leq \tau_7\| |\cdot|\|$ and $\| |\cdot|\| \leq \tau_8 \|\cdot\|$. Then we can obtain
 \begin{equation}
     ||x_k-x^*||\leq \|\psi^h_{k}\|\leq \tau_7 \||\psi^h_{k}|\| \leq \tau_7\left(\rho(P_{H3})+o(1)\right)^k\||\psi^h_{0}|\|.
 \end{equation}
 Hence, the convergence rate of DAGT-HB is only related to the spectral radius of $P_{H3}$. It is obvious $\rho(P_{H3})=\max \left\{\rho(P_{H4}),\rho\right\}$, where
 $$
P_{H4}=\left[\begin{array}{cc}
(1+\beta)I_N-\alpha C &-\beta I_N\\
I_N & 0
\end{array}\right].
 $$
 Based on \cite{ang2018heavy}, we can obtain that the optimal parameters are $\alpha=\frac{4}{(\sqrt{L_1}+\sqrt{\mu})^2}$ and $\beta=\frac{\sqrt{L_1}-\sqrt{\mu}}{\sqrt{L_1}+\sqrt{\mu}}$ and the spectral radius of $P_{H4}$ is $\frac{\sqrt{L_1}-\sqrt{\mu}}{\sqrt{L_1}+\sqrt{\mu}}$. Thus, for the quadratic optimization, the DAGT-HB can achieve a linear convergence rate with $\mathcal{O}\left(\max\left\{\rho^k,\left(\frac{\sqrt{L_1}-\sqrt{\mu}}{\sqrt{L_1}+\sqrt{\mu}}\right)^k\right\}\right)$.

  Finally, we establish the convergence rate of DAGT-NES for the quadratic optimization.  By invoking \eqref{nesx} and \eqref{nesy}, we can obtain
 { \begin{align}
      x_{k+1}&=x_k+\gamma(x_k-x_{k-1})-\alpha \left[\nabla_1 f\left(y_k, u_k\right)+\nabla \phi\left(y_k\right) s_k\right]\notag\\
      &=x_k+\gamma(x_k-x_{k-1})-\alpha [C(x_k+\gamma(x_k-x_{k-1}))+Hs_k]\notag\\
      &=(1+\gamma)(I_N-\alpha C)x_k-\gamma (I_N-\alpha C)x_{k-1}-\alpha H s_k.
  \end{align}}
 We   denote $\psi^n_{k}=\psi^h_{k}$, then we can derive
 \begin{align}
     \psi^n_{k+1}=P_{N}\psi^n_{k}+\theta_N,
 \end{align}
  where $P_N=$
$$
\left[\begin{array}{cccc}
(1+\gamma)(I_N-\alpha C) &-\gamma (I_N-\alpha C)&  0 & 0 \\
I_N & 0 & 0 & 0\\
 (1+2\gamma)(K_N-I_N)H&\gamma(I_N-K_N)H & A-K_N & 0\\
 0 & 0 &0 & A-K_N
\end{array}\right],
 $$
 and
 $$
\theta_N=\left[\begin{array}{c}
-\alpha H s_{k} \\
0\\
(1+\gamma)(H-{K}_NH)(x_{k+1}-x^*)\\
0
\end{array}\right] .
 $$
Similar to the above analysis, we have
\begin{align}
    P_{N1}\psi^n_{k+1} \leq P_{N2}\psi^n_{k},
    \label{PN1}
\end{align}
where 
 $$
P_{N1}=\left[\begin{array}{cccc}
I_N & 0& 0 & 0 \\
0 &I_N & 0 &0\\
(1+\gamma)( K_NH-H)&0& I_N & 0\\
 0& 0 & 0 & I_N
\end{array}\right],
 $$
 and $P_{N2}=$
$$
\left[\begin{array}{cccc}
(1+\gamma)(I_N-\alpha C) &-\gamma (I_N-\alpha C)&  0 & -\alpha H \\
I_N & 0 & 0 & 0\\
 (1+2\gamma)(K_N-I_N)H&\gamma(I-K_N)H & A-K_N & 0\\
 0 & 0 &0 & A-K_N
\end{array}\right].
 $$
 Further, we can derive
 \begin{align}
     \psi^n_{k+1} \leq P_{N1}^{-1}P_{N2}\psi^n_{k}=P_{N3}\psi^n_{k},
 \end{align}
where $P_{N3}=$
$$
\begin{aligned}
&\left[\begin{array}{cccc}
(1+\gamma)(I_N-\alpha C) &-\gamma (I_N-\alpha C)&  0 & -\alpha H \\
I_N & 0 & 0 & 0\\
 *&** & A-K_N & ***\\
 0 & 0 &0 & A-K_N
\end{array}\right].
\end{aligned}
 $$
 By employing Lemma 4, we can derive that there exists a matrix norm $\||\cdot|\|$ such that 
 \begin{align}
     \||\psi^n_{k+1}|\|\leq (\rho(P_{N3})+o(1)) \||\psi^n_{k}|\|.
 \end{align}
By  recalling that all norms are equivalent on finite-dimensional vector spaces, there always exist $\tau_9>0$ and $\tau_{10}>0$ such that $\|\cdot\| \leq \tau_9\| |\cdot|\|$ and $\| |\cdot|\| \leq \tau_{10} \|\cdot\|$. Then we can obtain
 \begin{equation}
     ||x_k-x^*||\leq \|\psi^n_{k}\|\leq \tau_9 \||\psi^n_{k}|\| \leq \tau_9\left(\rho(P_{N3})+o(1)\right)^k\||\psi^n_{0}|\|.
 \end{equation}
 Hence, the convergence rate of DAGT-NES is only related to the spectral radius of $P_{N3}$. It is obvious $\rho(P_{N3})=\max \left\{\rho(P_{N4}),\rho\right\}$, where
 $$
P_{N4}=\left[\begin{array}{cc}
(1+\gamma)(I_N-\alpha C) &-\gamma (I_N-\alpha C) \\
I_N & 0
\end{array}\right].
 $$
 Based on \cite{ang2018heavy}, we can obtain that the optimal parameters are $\alpha=\frac{4}{3L_1+\mu}$ and $\gamma=\frac{\sqrt{3\kappa+1}-2}{\sqrt{3\kappa+1}+2}$ and the spectral radius of $P_{N4}$ is $\frac{\sqrt{3\kappa+1}-2}{\sqrt{3\kappa+1}+2}$, where $\kappa=\frac{L_1}{\mu}$. Thus, for the quadratic optimization, the DAGT-NES can achieve a linear convergence rate with $\mathcal{O}\left(\max\left\{\rho^k,\left(\frac{\sqrt{3\kappa+1}-2}{\sqrt{3\kappa+1}+2}\right)^k\right\}\right)$. 

By summarizing the above analysis, we can obtain $\frac{\sqrt{3\kappa+1}-2}{\sqrt{3\kappa+1}+2}<\frac{\sqrt{L_1}-\sqrt{\mu}}{\sqrt{L_1}+\sqrt{\mu}}<\frac{L_1-\mu}{L_1+\mu}$. Thus, when $\rho\leq \frac{L_1-\mu}{L_1+\mu}$, the convergence rates of DAGT-HB and DAGT-NES are faster than DAGT.

Next, we further discuss the comparison among the convergence rate of algorithms DAGT-HB, DAGT-NES and  DAGT in the same case of $\alpha$. To facilitate the following analysis, we introduce the following lemma.

\begin{lemma}
\cite{ang2018heavy} If $\beta\geq (1-\sqrt{\alpha L_1})^2$, then $\rho(P_{H4})\leq \beta$. Additionally, if $\frac{1}{L_1}\leq \alpha \leq \frac{1}{\mu}$ and $\gamma\geq \frac{1-\sqrt{\alpha\mu }}{1+\sqrt{\alpha\mu}}$, then $\rho(P_{N4})\leq \sqrt{(1-\alpha\mu)\gamma}$.
\end{lemma}

We proceed by considering $\alpha=\frac{2}{\mu+L_1}$, which is the optimal step size for the DAGT algorithm. By selecting $\beta= (1-\sqrt{\alpha L_1})^2$ and $\gamma= \frac{1-\sqrt{\alpha\mu }}{1+\sqrt{\alpha\mu}}$, we can derive that $\rho(P_{H4})\leq \beta<\frac{L_1-\mu}{L_1+\mu}$ and $\rho(P_{N4})\leq \sqrt{(1-\alpha\mu)\gamma}<\frac{L_1-\mu}{L_1+\mu}$ based on the Lemma 16. Consequently, even when using the same $\alpha$, we can conclude that the convergence rates of the DAGT-HB and DAGT-NES algorithms are still faster than that of the DAGT algorithm.

Moreover, it is worth noting that the quadratic objective function optimization is a  well-established problem in trajectory planning for robots and autonomous driving \cite{Carnevale2022AggregativeFO,Carnevale2022Coordination,quirynen2024real}. Given its prevalence, it is highly relevant to delve deeper into our algorithm using this quadratic objective function as an example. Furthermore, for non-convex cost functions, our algorithms can be effectively combined with the Majorization Minimization (MM) algorithm framework \cite{sun2016majorization} to tackle the non-convex optimization problem. Specifically, we can find a surrogate function $g(x_k,x)$ that satisfies our assumptions and approximates the cost function at $x_k$. Then, we apply our algorithms to optimize the surrogate function $g(x_k,x)$, resulting in the next iteration point $x_{k+1}$. This process is repeated until the algorithms converge. The MM approach has been extensively utilized in various engineering applications \cite{figueiredo2007majorization,mairal2015incremental,oliveira2009adaptive}. We plan to conduct further theoretical analysis in our future work.

\begin{remark}
    In the quadratic optimization example, $h_i\geq 0$ is to ensure $\frac{\alpha}{N}H>0$ to obtain formulas \eqref{PD1}, \eqref{PH1} and \eqref{PN1}. For $h_i< 0$, this requires exploring other mathematical techniques to solve the problem, which will be considered in our future work. Besides, for the general distributed aggregrative optimization problem,  due to the shrinkage of some inequalities and coupling parameters, it is difficult to further compare the convergence rate only through simulation results\cite{wang2024momentum}\cite{chen2024achieving}. 
\end{remark}

\section{Numerical Simulation}
In this section, we  perform the following two simulations to verify the effectiveness of our proposed methods  DAGT-HB and DAGT-NES.
\vspace{-1 em}
\subsection{Simulation 1: The Optimal Placement Problem}

In an optimal placement problem in $\mathbb{R}^2$, suppose that there are 5 entities which are located at $r_1=(10,4)$, $r_2=(1,3)$, $r_3=(2,7)$, $r_4=(8,10)$ and $r_5=(3,9)$. And there are 5 free entities, each of which can privately be accessible to some of the fixed 5 entities. The purpose is to determine the optimal position $x_i$, $i\in \{1,2,3,4,5\}$  of the free entity so as to minimize the sum of all distances from the current position of each free entity to the corresponding fixed entity location and the distances from each entity to the weighted center of all free entities. Therefore, the cost function of  each free entity can be  modeled as follows:
\begin{equation}
    f_i(x_i,u(x))=\omega_i ||x_i-r_i||^2+ ||x_i-u(x)||^2, \quad i=1,\ldots,5,
\end{equation}
where $\omega_i$ represents the weight and is set to $20$.  We set $u(x)=\frac{\sum_{i=1}^5x_i}{5}$. So in this condition, $\phi_i$ is the identity mapping for $i=1,\ldots,5$. The communication
graph is randomly chosen to be  connected and doubly stochastic. 

We select the initial point $x_{1,-1}=(0,11)$, $x_{2,-1}=(9,8)$, $x_{3,-1}=(9,1)$, $x_{4,-1}=(1,4)$ and $x_{5,-1}=(3,1)$; $x_{1,0}=y_{1,0}=(2,9)$, $x_{2,0}=y_{2,0}=(8,6)$, $x_{3,0}=y_{3,0}=(7,3)$, $x_{4,0}=y_{4,0}=(4,7)$ and $x_{5,0}=y_{5,0}=(8,3)$. We initialize $u_{i,0}=\phi_i(x_{i,0})$ and $s_{i,0}=\nabla_2 f_i(x_{i,0},u_{i,0})$ for $i=1,2,\ldots,5$. {Based on Theorem 1 and Theorem 2, we derive the upper bounds $\alpha<0.00639, \beta<0.00943$, and $\alpha<0.00572, \gamma<0.00817$. Subsequently, we elect to set $\alpha=0.005, \beta=0.009$, and $\gamma=0.008$ respectively.}

Then we use the DAGT-HB method to solve the optimal placement problem and the results are shown in Fig. \ref{fig1}. In Fig. \ref{fig1} (a)-(b),  we can see that all agents can converge to their best positions and the optimal positions are $x_1^*=(9.7524,4.1248)$, $x_2^*=(1.810,3.1714)$, $x_3^*=(2.1333,6.9810)$, $x_4^*=(7.8416,9.8381)$ and $x_5^*=(3.0857,8.8857)$ respectively.
Fig. \ref{fig1} (c)-(d)  show the  evolution of $u_{i,k}$, indicating that  the estimate $u_{i,k}$ of each free entity converges to optimal aggregative position $u(x^*)=(4.8,6.6)$ with a rapid speed.  Furthermore, Fig. \ref{fig1} implies that the convergence rate is fast, which supports our theoretical analysis.  Meanwhile, the evolutions of $x_{i,k}$ and $u_{i,k}$ for DAGT-NES are similar to those for DAGT-HB and the results are shown in Fig. \ref{fig2}.  

In order to demonstrate the superiority of DAGT-HB and DAGT-NES, we compare them with DAGT and A-DAGT \cite{chen2024achieving} under the same initial conditions. The residual error is shown in Fig. \ref{fig3}. We find that the convergence speeds of DAGT-HB and DAGT-NES are significantly faster than that of DAGT, indicating that the introduced momentum term $x_{i,k}-x_{i,k-1}$ can enhance the convergence rate of the algorithm. Furthermore, compared with A-DAGT, the convergence rates of our algorithms are faster. These results demonstrate the superiority of our algorithms.


{

Additionally, we conduct the simulation of DAGT-HB and DAGT-NES with various parameters: specifically, $\beta=0$ (DAGT), $\beta=0.005, \beta=0.01,$ and $\beta=0.05$, and $\gamma=0.005, \gamma=0.01,$ and $\gamma=0.05$ respectively, with $\alpha=0.005$. The initial value is consistent with the above simulations, and the relative error trajectories of $(F(x_k)-F(x^*))^2$ are illustrated in Fig. \ref{fig4}. In Fig. \ref{fig4}, it is evident that the convergence rates of both DAGT-HB and DAGT-NES accelerate with the increase of $\beta$ and $\gamma$. This observation implies that our algorithms offer a wider range of parameters to choose from, indicating a greater potential to accelerate the convergence rate. 

}

\begin{figure}
    \centering
    \includegraphics[width=\linewidth]{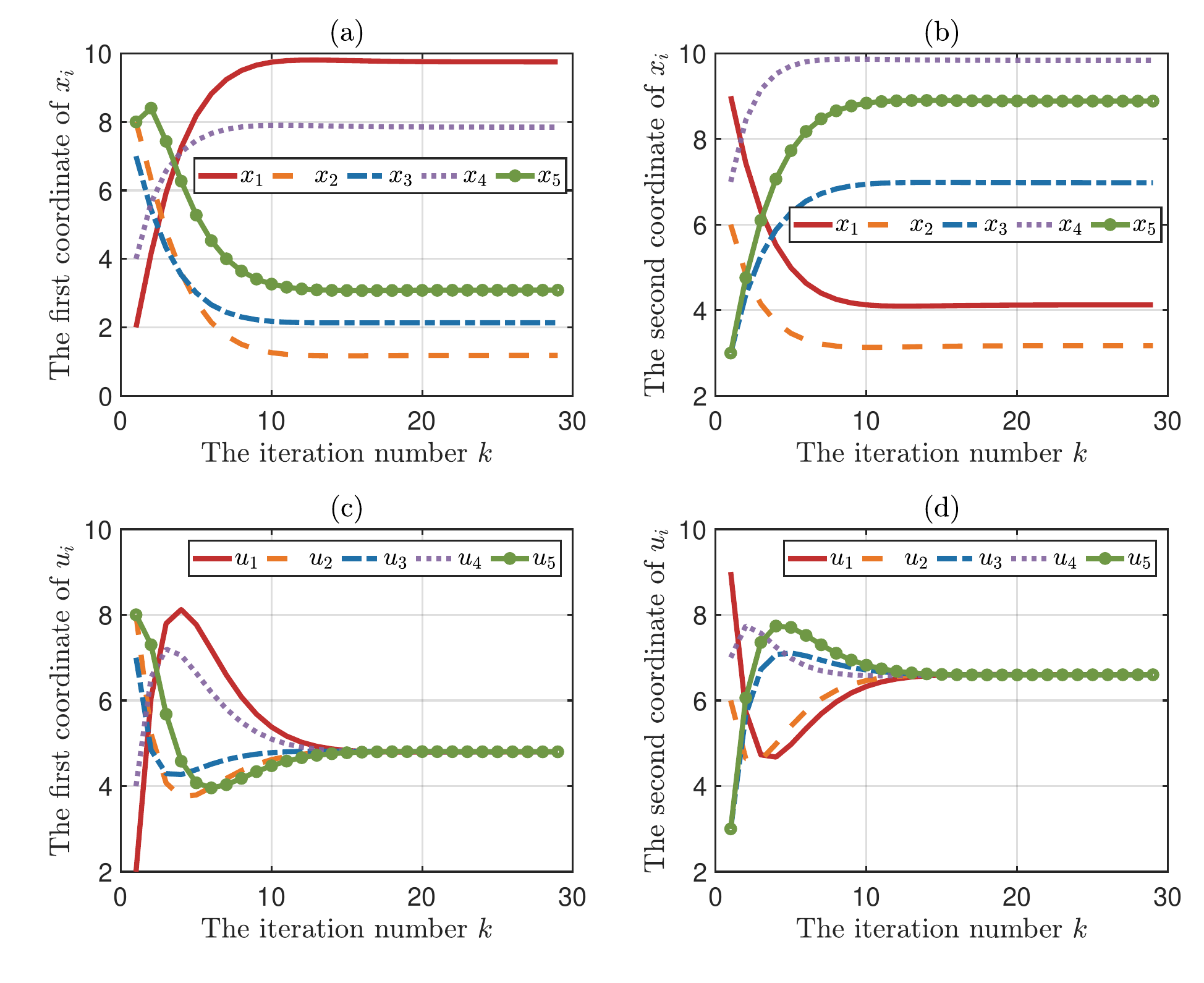}
    \caption{(a):  The evolution of the first coordinate of $x_{i,k}$ by using DAGT-HB. (b): The evolution of the second coordinate of $x_{i,k}$ by using DAGT-HB. (c): The evolution of the first coordinate of $u_{i,k}$ by using DAGT-HB. (d): The evolution of the second coordinate of $u_{i,k}$ by using DAGT-HB.}
    \label{fig1}
\end{figure}

\begin{figure}
    \centering
    \includegraphics[width=\linewidth]{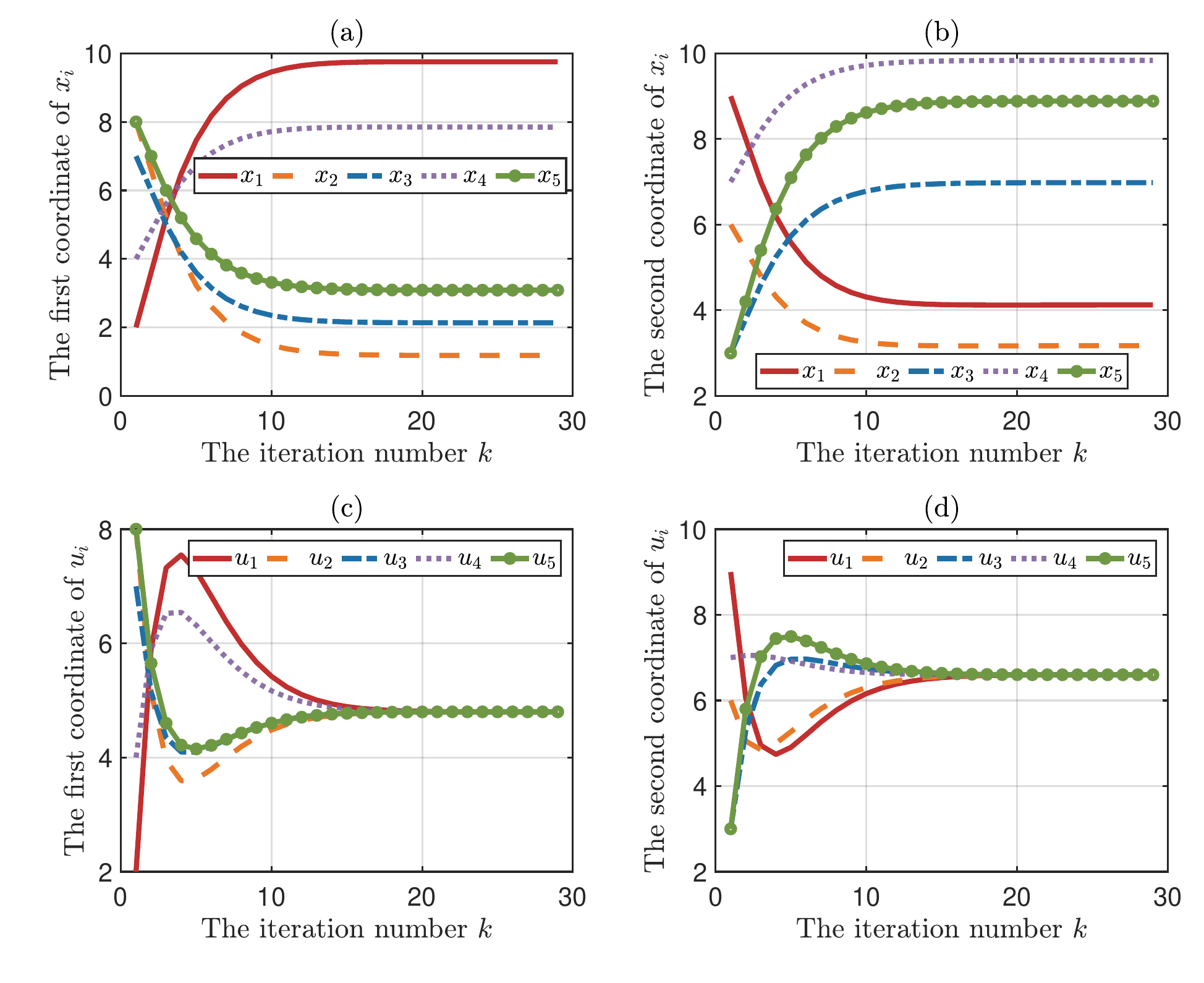}
    \caption{(a):  The evolution of the first coordinate of $x_{i,k}$ by using DAGT-NES. (b): The evolution of the second coordinate of $x_{i,k}$ by using DAGT-NES. (c): The evolution of the first coordinate of $u_{i,k}$ by using DAGT-NES. (d): The evolution of the second coordinate of $u_{i,k}$ by using DAGT-NES.}
    \label{fig2}
\end{figure}



\begin{figure}
    \centering
    \includegraphics[width=\linewidth]{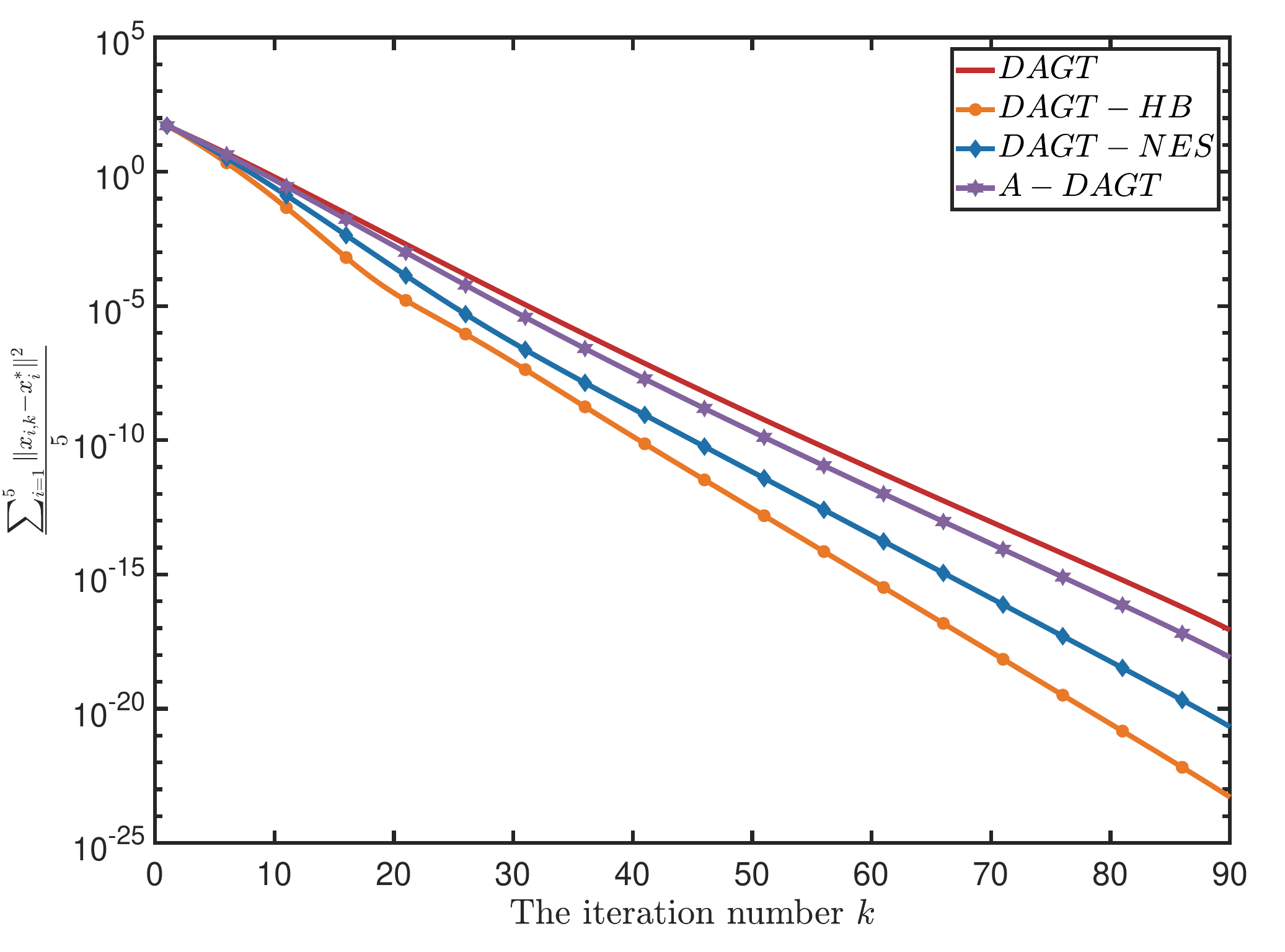}
    \caption{The state error $\frac{\sum_{i=1}^5\|x_{i,k}-x_i^*\|^2}{5}$ comparison among DAGT, DAGT-HB, DAGT-NES and A-DAGT.}
    \label{fig3}
\end{figure}


\begin{figure}
    \centering
    \includegraphics[width=\linewidth]{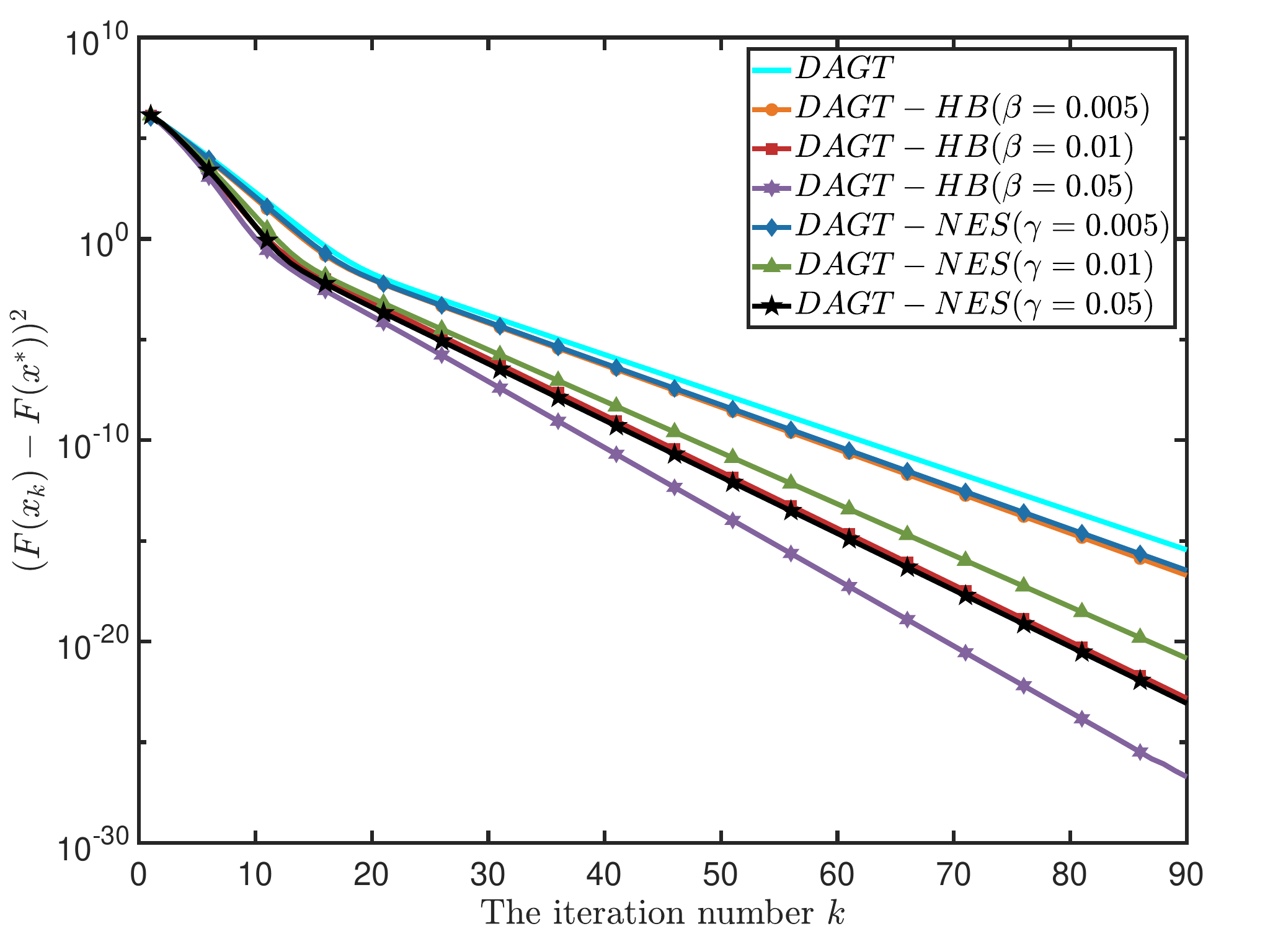}
    \caption{The relative error trajectories of $(F(x_k)-F(x^*))^2$ under different parameters.}
    \label{fig4}
\end{figure}
\vspace{-1 em}
\subsection{Simulation 2: A Class of Nash-Cournot Game of Generation Systems \cite{chen2024achieving}}
Here, we investigate the economic dispatch problem of generation systems in smart grids, where one of the most conspicuous features is the dependence of electricity prices on the overall output power demand. In the standard Nash-Cournot game, each generation system is required to determine the output power value to minimize its individual cost function. However, for the entire power market, each system needs to collaboratively optimize the overall social cost, which is the distributed aggregative optimization problem. We investigate the DAGT-HB and DAGT-NES to solve the class of Nash-Cournot game played by $N=50$ generation systems. For each generation system $i=1,\ldots,50$, the local cost function can be formulated as follows:
\begin{equation}
    f_i\left(x_i,u(x)\right)=\kappa_i x_i^2+\theta_i x_i+\sigma_i-(\omega_1-\omega_2 u(x))x_i,
    \label{sim2}
\end{equation}
where $\kappa_i, \theta_i, \sigma_i$ are characteristic parameters of generation system  $i$, $u(x)=\sum_{i=1}^N x_i$ and $\omega_1, \omega_2$ are constants.  Therefore, without considering mechanical
and electromagnetic losses, the goal of this work is that each
generation system collaborates to obtain an optimal solution
that makes the social cost optimal.

The generation system $i$’s characteristics are selected randomly
at the following intervals: $\kappa_i \in [0.5,2.5], \theta_i \in [10,20]$ and $\sigma_i \in [5,20]$. We set randomly the initial point $x_{i,-1}=x_{i,0}$ at the following intervals: $x_{i,0} \in [50,100]$ for $i=1,2,\ldots,50$. Further, we initialize $u_{i,0}=\phi_i(x_{i,0})$ and $s_{i,0}=\nabla_2 f_i(x_{i,0},u_{i,0})$ for $i=1,2,\ldots,50$ and $\omega_1=200$ and $\omega_2=0.01$. Based on Theorem 1 and Theorem 2, we derive the upper bounds $\alpha<0.00347, \beta<0.00643$, and $\alpha<0.00312, \gamma<0.00528$. Then we set  the step size $\alpha=0.003$ and choose the momentum term $\beta=0.006$ and $\gamma=0.005$ respectively.

We employ the DAGT-HB and DAGT-NES algorithms to solve the problem \eqref{sim2}. We take $x_{1,k},x_{11,k},x_{21,k},x_{31,k}$ and $x_{41,k}$ as examples and the evolutions of $x_{i}s'$ by leveraging the DAGT-HB and DAGT-NES algorithms are shown in Fig. \ref{fig5}. Fig. \ref{fig5} shows that $x_{i,k}$ rapidly converges to the same optimal value by using the DAGT-HB and DAGT-NES algorithms respectively, which supports our theoretical analysis.  Furthermore,   we compare our algorithms with DAGT and A-DAGT under the same initial conditions.  The residual error $\frac{\sum_{i=1}^{50}\|x_{i,k}-x_i^*\|^2}{50}$ is depicted in Fig. \ref{fig6}. We find that the convergence speeds of DAGT-HB and DAGT-NES are significantly faster than that of DAGT and A-DAGT, which demonstrates the superiority of DAGT-HB and DAGT-NES.

Subsequently, we investigate the effects of adjusting the momentum parameters $\beta$ and $\gamma$ on the convergence rates of the DAGT-HB and DAGT-NES algorithms, respectively. With identical initial conditions, we examine the number of iterations required for both algorithms to converge to $10^{-6}$ for various values of $\beta$ and $\gamma$. The results are illustrated in Fig. \ref{fig7}. As shown in Fig. \ref{fig7}, the convergence speeds of both algorithms initially increase and then decrease as $\beta$ and $\gamma$ increase. Specifically, the DAGT-HB algorithm achieves its maximum convergence rate when $\beta$ lies between 0.4 and 0.5, while the DAGT-NES algorithm reaches its maximum convergence rate when $\gamma$ falls between 0.5 and 0.6.

Furthermore, we explore the impact of different topology graph structures on the performance of the DAGT-HB and DAGT-NES algorithms. Three distinct topological graph structures are considered: star, ring, and fully connected \cite{zhu2022topology}. We set the step size $\alpha=0.003$ and choose the momentum terms $\beta=0.006$ and $\gamma=0.005$. The results are presented in Fig. \ref{fig8}. As evident from Fig. \ref{fig8}, increased interaction among agents accelerates the convergence of both algorithms. However, this improvement comes at the cost of additional communication overhead, highlighting a trade-off.

To assess the robustness of our algorithms, we evaluate their performance under scenarios with communication delay and communication error. In the communication delay scenario, we introduce a two-step delay for every communication among agents. The results are shown in the left subgraph of Fig. \ref{fig9}. Both DAGT-HB and DAGT-NES algorithms demonstrate a relatively fast convergence rate, indicating their robustness to communication delays. Moreover, compared to A-DAGT and DAGT, our algorithms converge more quickly. In the communication error scenario, we introduce white noise interference in every communication for the DAGT-HB algorithm. The results are depicted in the right subgraph of Fig. \ref{fig9}. Despite the presence of noise, both algorithms still exhibit a relatively fast convergence rate, further supporting their robustness. Additionally, in the same setting, we conducte the comparative experiment and find that our algorithms converge faster than DAGT and A-DAGT algorithms. Overall, these results provide strong evidence for the effectiveness and robustness of our proposed algorithms.

\begin{figure}
    \centering
    \includegraphics[width=\linewidth]{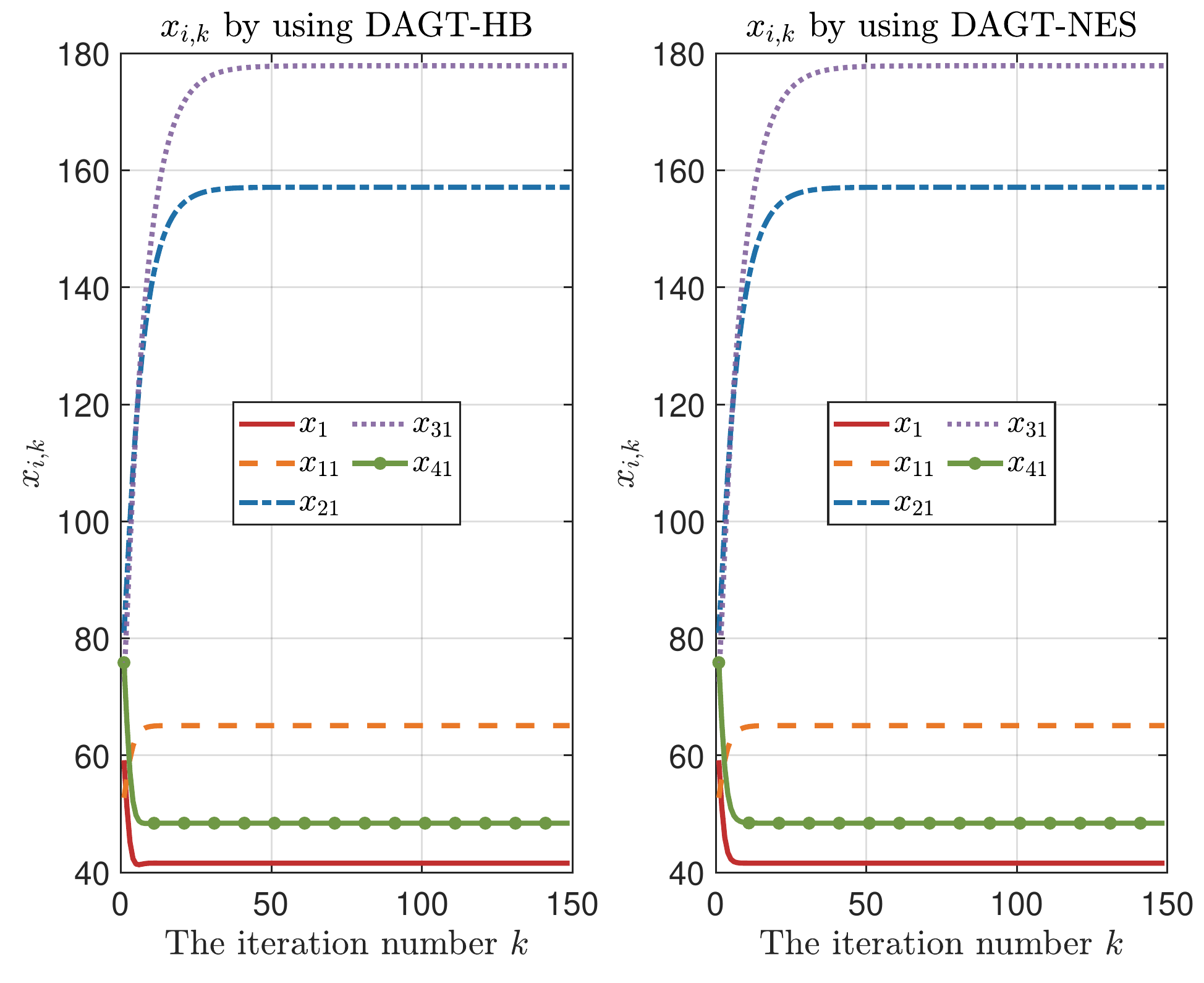}
    \caption{Left: The evolution of $x_{i,k}$ by using DAGT-HB. Right: The evolution of $x_{i,k}$ by using DAGT-NES.}
    \label{fig5}
\end{figure}

\begin{figure}
    \centering
    \includegraphics[width=\linewidth]{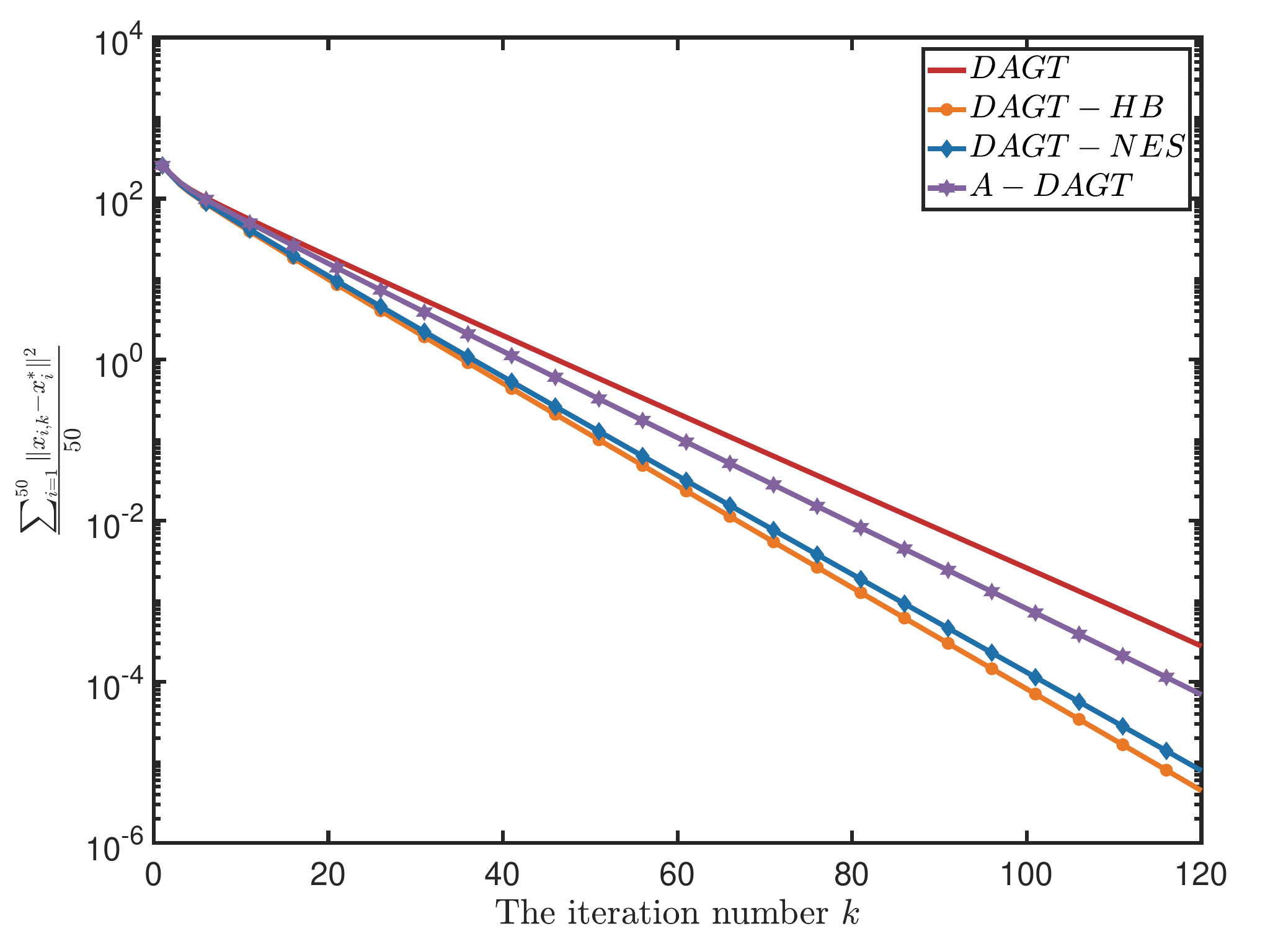}
    \caption{The state error  $\frac{\sum_{i=1}^{50}\|x_{i,k}-x_i^*\|^2}{50}$ comparison among DAGT, DAGT-HB, DAGT-NES and A-DAGT.}
    \label{fig6}
\end{figure}

\begin{figure}
    \centering
    \includegraphics[width=\linewidth]{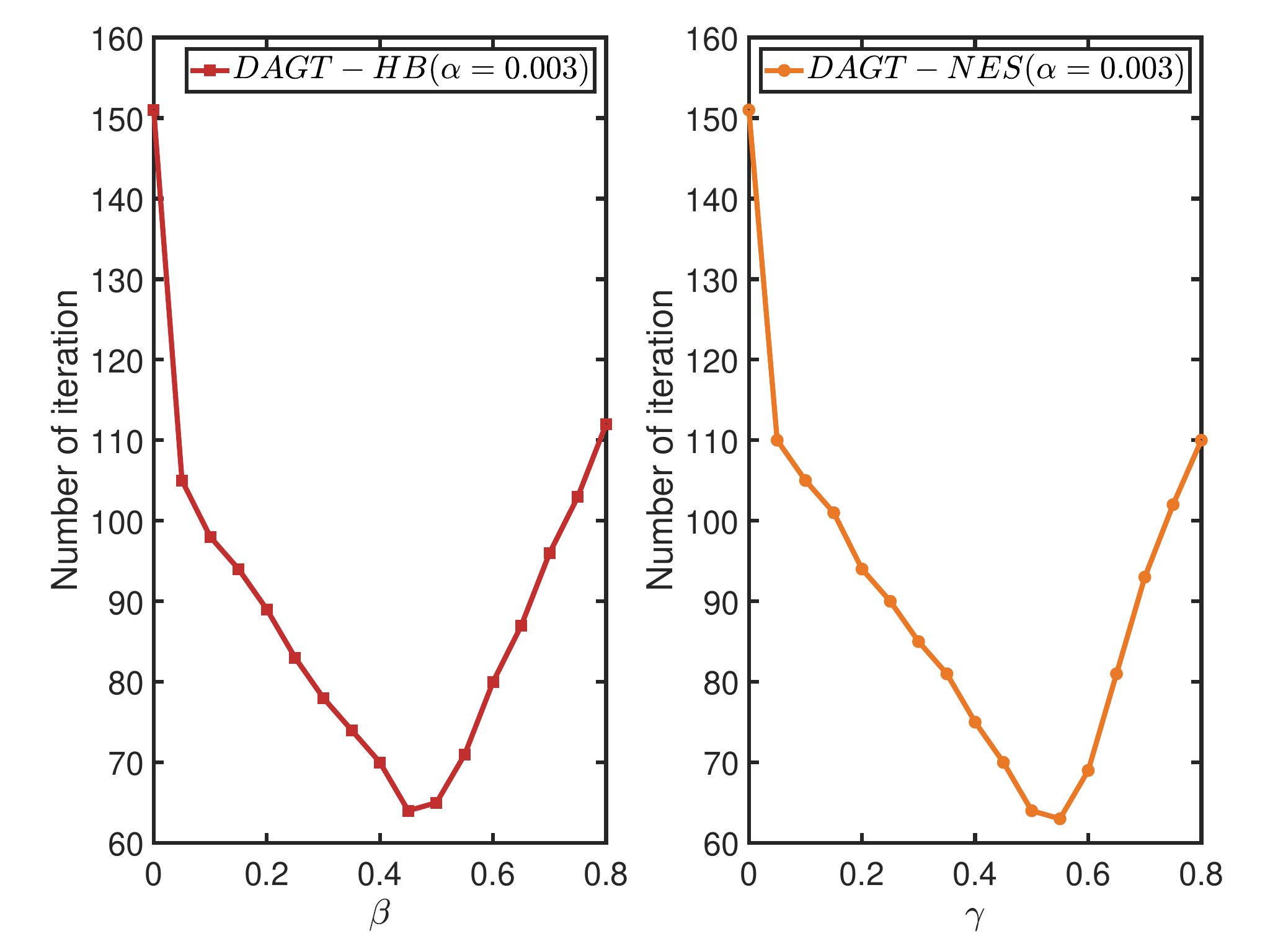}
    \caption{Effect of adjusting parameters $\beta$ and $\gamma$ for DAGT-HB and DAGT-NES, respectively.}
    \label{fig7}
\end{figure}

\begin{figure}
    \centering
    \includegraphics[width=\linewidth]{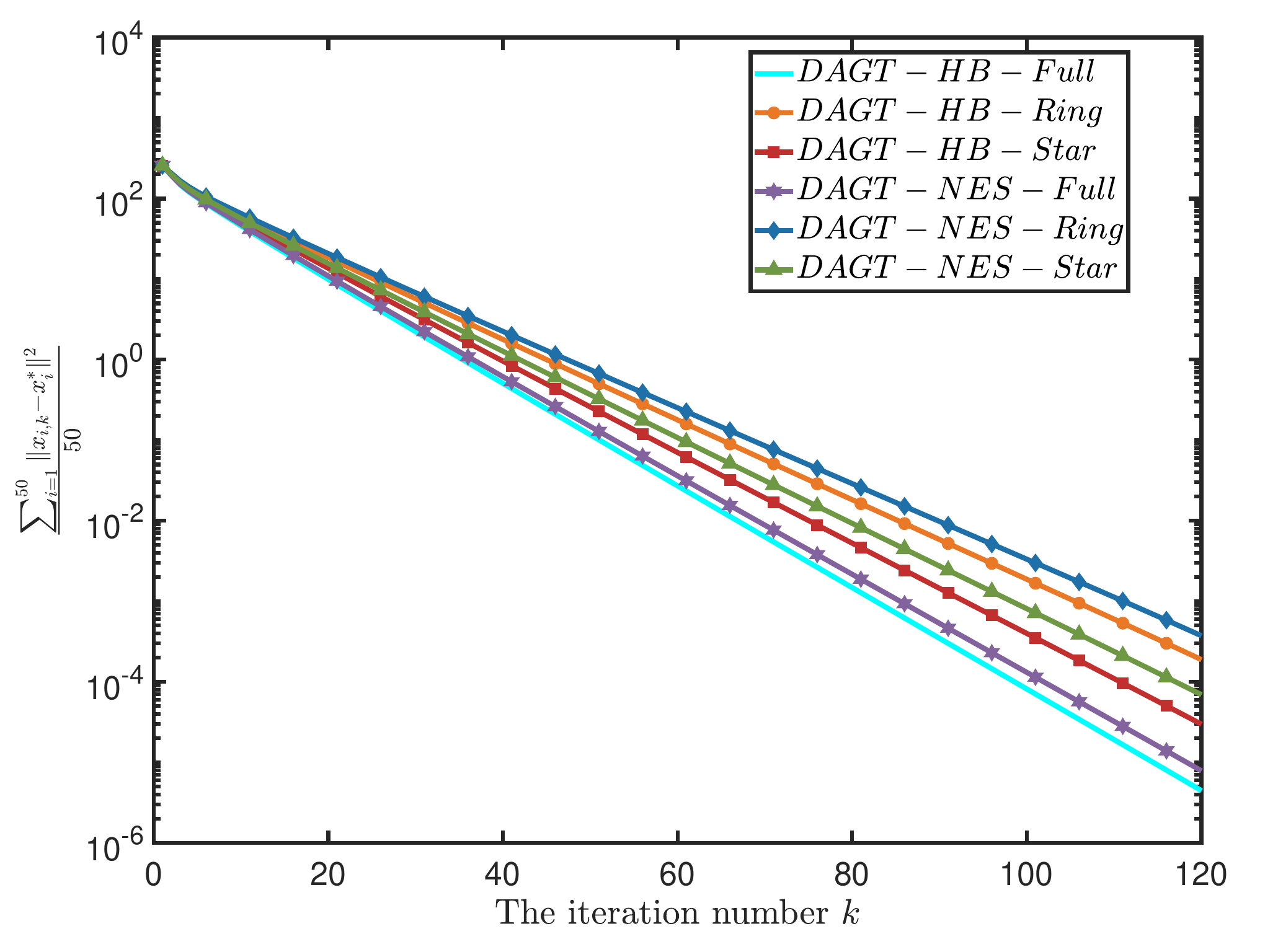}
    \caption{The residual  error  comparison over different graphs.}
    \label{fig8}
\end{figure}

\begin{figure}
    \centering
    \includegraphics[width=\linewidth]{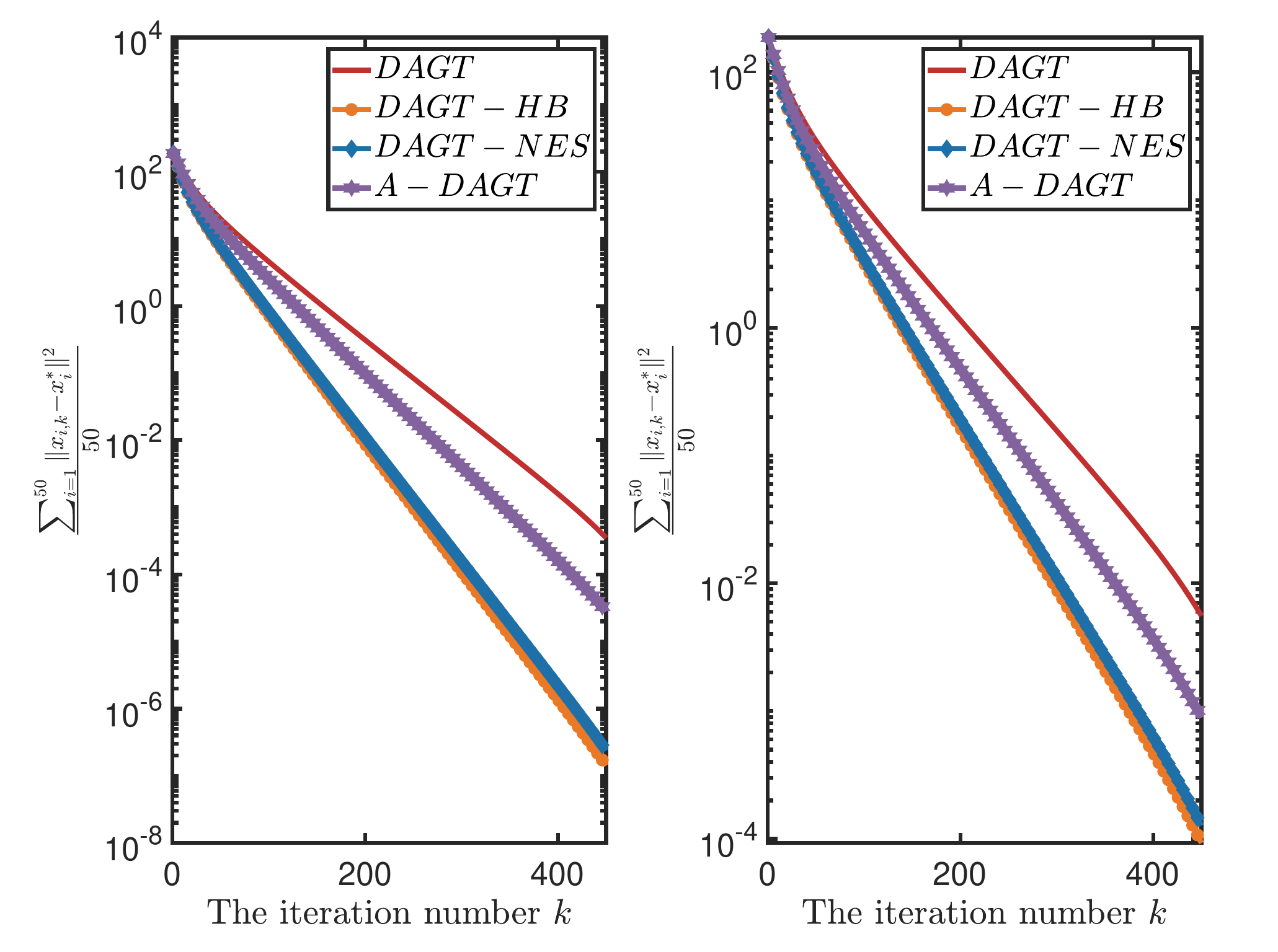}
    \caption{Left: The residual  error  comparison with two step  communication delay. Right: The residual  error  comparison with communication error.}
    \label{fig9}
\end{figure}

\vspace{-1 em}
\section{Conclusion}

This paper proposes two novel algorithms called DAGT-HB and DAGT-NES to solve the distributed aggregative
optimization problem in a network.   Inspired by the accelerated algorithms, we combine heavy ball and Nesterov’s accelerated method with the distributed
aggregative gradient tracking method. Furthermore, we show that the algorithms DAGT-HB and DAGT-NES can converge
to an optimal solution at a global linear convergence rate when
the objective function is smooth and strongly convex and when the step size and momentum term can be selected appropriately. Additionally,  within the context of the quadratic distributed aggregative optimization problem, we illustrate that the convergence rates of DAGT-HB and DAGT-NES surpass that of DAGT when the spectral radius of the weighted adjacency matrix falls within a specific range. Finally, we use two numerical simulations
to verify the effectiveness and
superiority of DAGT-HB and DAGT-NES. Under the same conditions, DAGT-HB and DAGT-NES can achieve much faster convergence of the state error and cost function compared to the vanilla DAGT method and A-DAGT algorithm.  Future study may focus on  the sensitivity of parameters selection in DAGT-HB and DAGT-NES and the extension of DAGT-HB and DAGT-NES  to unbalanced graph, non-convex objective function and constrained distributed aggregation optimization problems.
\vspace{-1 em}
\section{APPENDIX}

\subsection{Proof of Lemma 11}
Denote the equilibrium point of \eqref{nesx}-\eqref{ness} as $x^*=$col $\left(x_1^*, \ldots, x_N^*\right)$, $y^*=$col $\left(y_1^*, \ldots, y_N^*\right)$, $u^*=$ col$\left(u_1^*, \ldots, u_N^*\right)$, and $s^*=\operatorname{col}\left(s_1^*, \ldots, s_N^*\right)$. By \eqref{nesy}, we know $y^*=x^*$. Next, from \eqref{nesx}, \eqref{nesu} and \eqref{ness}, we can derive
\begin{align}
    & \nabla_1 f\left(y^*, u^*\right)+\nabla \phi\left(y^*\right) s^*=\mathbf{0}_{Nd}, \label{lem11 1}\\
& \mathcal{L} u^*=\mathbf{0}_{Nd},\quad \mathcal{L} s^*=\mathbf{0}_{Nd} ,
\end{align}
where $\mathcal{L}=L\otimes I_d$ and $L$ is Laplacian matrix of graph $\mathcal{G}$. Because of the property of the Laplace matrix, so $u_i^*$ must be equal to $u_j^*$ and $s_i^*$ must be equal to $s_j^*$ for all $i\neq j$. Due to equality  \eqref{mean u} and \eqref{mean s}, it leads to
\begin{align}
    u_i^* & =\frac{1}{N} \sum_{i=1}^N \phi_i\left(y_{i}^*\right)=u(y^*), \label{lem11 2}\\
s_i^* & =\frac{1}{N} \sum_{i=1}^N \nabla_2 f_i\left(y_{i}^*, u(y^*)\right)  .\label{lem11 3}
\end{align}
By inserting \eqref{lem11 2} and \eqref{lem11 3} into \eqref{lem11 1}, we can obtain $\nabla F(x^*)=\nabla F(y^*)=0$. Because $F$ is $\mu$-strongly convex, $x^*$ is the unique optimal solution to the problem (1).

\vspace{-1 em}

\subsection{Proof of Lemma 12}
For $\left\|x_{k+1}-x^*\right\|$, by invoking \eqref{nesx}, it leads to
\begin{equation}
    \begin{aligned}
        & \left\|x_{k+1}-x^*\right\| \\
& =\left\|y_k-y^*-\alpha\left[\nabla_1 f\left(y_k, u_k\right)+\nabla \phi\left(y_k\right) s_k\right]\right\| \\
& \leq \Big\Vert y_k-y^*-\alpha\left[\nabla_1 f\left(y_k, \mathbf{1}_N \otimes \bar{u}_k\right)\right. \\
&\quad \left.+\nabla \phi\left(y_k\right) \mathbf{1}_N \otimes \frac{1}{N} \sum_{i=1}^N \nabla_2 f_i\left(y_{i, k},  \bar{u}_k\right)\right]+\alpha \nabla f\left(y^*\right) \Big\Vert \\
&\quad +\alpha \| \nabla_1 f\left(y_k, u_k\right)+\nabla \phi\left(y_k\right) \mathbf{1}_N \otimes \bar{s}_k-\nabla_1 f\left(y_k, \mathbf{1}_N \otimes \bar{u}_k\right) \\
&\quad -\nabla \phi\left(y_k\right) \mathbf{1}_N \otimes \frac{1}{N} \sum_{i=1}^N \nabla_2 f_i\left(y_{i, k},  \bar{u}_k\right) \| \\
&\quad +\alpha\left\|\nabla \phi\left(y_k\right) s_k-\nabla \phi\left(y_k\right) \mathbf{1}_N \otimes \bar{s}_k\right\|. \\
    \end{aligned}
    \label{lem12 1}
\end{equation}
From Lemma $1$, we can bound the first term of the right term of \eqref{lem12 1} as follows:
\begin{equation}
    \begin{aligned}
        &  \Big\Vert y_k-y^*-\alpha\left[\nabla_1 f\left(y_k, \mathbf{1}_N \otimes \bar{u}_k\right)\right. \\
& \quad\left.+\nabla \phi\left(y_k\right) \mathbf{1}_N \otimes \frac{1}{N} \sum_{i=1}^N \nabla_2 f_i\left(y_{i, k},  \bar{u}_k\right)\right]+\alpha \nabla F\left(y^*\right) \Big\Vert \\
&  \leq(1-\mu \alpha)\left\|y_k-y^*\right\|.
    \end{aligned}
    \label{lem12 2}
\end{equation}
Then by inserting \eqref{nesy} into \eqref{lem12 2} and noting $x^*=y^*$  we can obtain
\begin{equation}
 \begin{aligned}
     \left\|y_k-y^*\right\|\leq \left\|x_k-x^*\right\|+ \gamma\left\|x_k-x_{k-1}\right\|.
 \end{aligned}   
\end{equation}
For the second term, since $\nabla_1 f(x, {y})+\nabla \phi(x) \mathbf{1}_{{N}} \otimes \frac{1}{{N}} \sum_{i=1}^{{N}} \nabla_2 f_i\left(x_i, y_i\right)$ is  $L_1$-Lipschitz continuous  according to Assumption 2 and $\mathbf{1}_N \otimes \bar{u}_k=\mathcal{K} u_k$,  we can get
\begin{equation}
    \begin{aligned}
        & \alpha \| \nabla_1 f\left(y_k, u_k\right)+\nabla \phi\left(y_k\right) \mathbf{1}_N \otimes \bar{s}_k-\nabla_1 f\left(y_k, \mathbf{1}_N \otimes \bar{u}_k\right) \\
&\quad -\nabla \phi\left(y_k\right) \mathbf{1}_N \otimes \frac{1}{N} \sum_{i=1}^N \nabla_2 f_i\left(y_{i, k},  \bar{u}_k\right) \| \\
& \leq \alpha L_1\left\|u_k-\mathcal{K} u_k\right\|.\\
    \end{aligned}
\end{equation}
For the last term, by using Assumption $4$ and $\mathbf{1}_N \otimes \bar{s}_k=\mathcal{K} s_k$ we can  the following inequality:
\begin{equation}
    \alpha\left\|\nabla \phi\left(y_k\right) s_k-\nabla \phi\left(y_k\right) \mathbf{1}_N \otimes \bar{s}_k\right\| \leq \alpha L_3\left\|s_k-\mathcal{K} s_k\right\|.
    \label{lem12 3}
\end{equation}
Then by bonding \eqref{lem12 2}-\eqref{lem12 3} with \eqref{lem12 1}, we complete the proof.

\vspace{-1 em}
\subsection{Proof of Lemma 13}
For $\left\|x_{k+1}-x_k\right\|$, by invoking \eqref{nesx} and noting  $\nabla F(x^*)=0$,  we have
    \begin{align}
        & \left\|x_{k+1}-x_k\right\| \notag\\
& =\left\|\gamma(x_k-x_{k-1})-\alpha(\nabla_1 f\left(y_k, u_k\right)+\nabla \phi\left(y_k\right) s_k)\right\| \notag\\
& \leq \alpha \Big\Vert \nabla_1 f\left(y_k, u_k\right)+\nabla \phi\left(y_k\right) \mathcal{K} s_k-\nabla_1 f\left(y^*, \mathbf{1}_N \otimes \bar{u}^*\right) \notag\\
& \quad\quad-\nabla \phi\left(x^*\right)\left[\mathbf{1}_N \otimes \frac{1}{N} \sum_{i=1}^N \nabla_2 f_i\left(y^*_i,  \bar{u}^*\right)\right] \Big\Vert \notag\\
& \quad+\alpha\left\|\nabla \phi\left(y_k\right)\left(s_k-\mathcal{K} s_k\right)\right\|+\gamma\left\|x_k-x_{k-1}\right\|. 
\label{lem13 1}
    \end{align}
By utilizing Assumption $2$ and triangle inequality of norm, we can obtain the following formula:
    \begin{align}
        &   \Big\Vert \nabla_1 f\left(y_k, u_k\right)+\nabla \phi\left(y_k\right) \mathcal{K} s_k-\nabla_1 f\left(y^*, \mathbf{1}_N \otimes \bar{u}^*\right) \notag\\
& \quad-\nabla \phi\left(y^*\right)\left[\mathbf{1}_N \otimes \frac{1}{N} \sum_{i=1}^N \nabla_2 f_i\left(y^*_i,  \bar{u}^*\right)\right] \Big\Vert \notag\\
& \leq  L_1\left(\left\|y_k-x^*\right\|+\left\|u_k-\mathbf{1}_N \otimes \bar{u}^*\right\|\right)\notag\\
& \leq  L_1\left(\left\|x_k-x^*\right\|+\left\|u_k-\mathcal{K} u_k\right\|\right) \notag\\
& \quad+\gamma L_1\left\|x_k-x_{k-1}\right\|+L_1\left\|\mathcal{K} u_k-\mathbf{1}_N \otimes \bar{u}^*\right\|.
\label{lem13 2}
    \end{align}
By \eqref{8 big2}, we know
\begin{equation}
    \begin{aligned}
        \left\|\mathcal{K} u_k-\mathbf{1}_N \otimes \bar{u}^*\right\| & \leq L_3\left\|y_k-y^*\right\|\\
        & \leq  L_3\left\|x_k-x^*\right\|+ L_3\gamma\left\|x_k-x_{k-1}\right\|.
    \end{aligned}
\end{equation}
Then by using Assumption $4$ we can obtain
\begin{equation}
    \left\|\nabla \phi\left(y_k\right)\left(s_k-\mathcal{K} s_k\right)\right\| \leq L_3 \left\|s_k-\mathcal{K} s_k\right\|.
    \label{lem13 3}
\end{equation}
Finally, by inserting \eqref{lem13 2}-\eqref{lem13 3}  into \eqref{lem13 1}, we can obtain the Lemma $13$.
\vspace{-1 em}
\subsection{Proof of Lemma 14}
For $\left\|u_{k+1}-\mathcal{K} u_{k+1}\right\|$, by invoking \eqref{nesu}, it leads to
\begin{equation}
\begin{aligned}
& \left\|u_{k+1}-\mathcal{K} u_{k+1}\right\| \\
& =\left\|\mathcal{A} u_k+\phi\left(y_{k+1}\right)-\phi\left(y_k\right)-\mathcal{K} \mathcal{A} u_k-\mathcal{K}\left[\phi\left(y_{k+1}\right)-\phi\left(y_k\right)\right]\right\| \\
& \leq \rho\left\|u_k-\mathcal{K} u_k\right\|+\|I_{Nd}-\mathcal{K}\|\left\|\phi\left(y_{k+1}\right)-\phi\left(y_k\right)\right\| \\
& \leq \rho\left\|u_k-\mathcal{K} u_k\right\|+L_3\left\|y_{k+1}-y_k\right\|\\
& = \rho\left\|u_k-\mathcal{K} u_k\right\|+L_3\left\|(1+\gamma)(x_{k+1}-x_k)-\gamma (x_k-x_{k-1})\right\|\\
& \leq\rho\left\|u_k-\mathcal{K} u_k\right\|+L_3(1+\gamma)\left\|x_{k+1}-x_k\right\|+\gamma L_3\left\|x_k-x_{k-1}\right\|,\\
\end{aligned}
\label{lem14 1}
\end{equation}
where Lemma $2$ has been utilized to obtain the first inequality, and by using  Assumption $4$ we can obtain the third inequality. Then by substituting \eqref{lem14} into \eqref{lem14 1} we can obtain the Lemma $14$.
\vspace{-1 em}
\subsection{Proof of Lemma 15}
 For $\left\|s_{k+1}-\mathcal{K} s_{k+1}\right\|$, by invoking \eqref{ness}  we can obtain
\begin{equation}
\begin{aligned}
& \left\|s_{k+1}-\mathcal{K} s_{k+1}\right\| \\
& =||\mathcal{A} s_k+\nabla_2 f\left(y_{k+1}, u_{k+1}\right)-\nabla_2 f\left(y_k, u_k\right)\\
&\quad-\mathcal{K}\mathcal{A} s_k-\mathcal{K}[\nabla_2 f\left(y_{k+1}, u_{k+1}\right)-\nabla_2 f\left(y_k, u_k\right)]|| \\
& \leq \left\|\mathcal{A}s_k-\mathcal{K} s_k\right\|\\
& \quad+\|I_{Nd}-\mathcal{K}\|\left\|\nabla_2 f\left(y_{k+1}, u_{k+1}\right)-\nabla_2 f\left(y_k, u_k\right)\right\|\\
& \leq \rho\left\|s_k-\mathcal{K} s_k\right\|+\left\|\nabla_2 f\left(y_{k+1}, u_{k+1}\right)-\nabla_2 f\left(y_k, u_k\right)\right\| \\
& \leq \rho\left\|s_k-\mathcal{K} s_k\right\|+  L_2\left(\left\|y_{k+1}-y_k\right\|+\left\|u_{k+1}-u_k\right\|\right),\\
\end{aligned}
\label{lem15 1}
\end{equation}
where Assumption $3$ has been leveraged in the last inequality. Then by \eqref{10 big1}, we know
\begin{equation}
\begin{aligned}
    & \left\|u_{k+1}-u_k\right\|
 \leq 2 \left\|u_k-\mathcal{K} u_k\right\|+ L_3\left\|y_{k+1}-y_k\right\|.
    \end{aligned}
    \label{lem15 2}
\end{equation}
By substituting  \eqref{lem15 2} into \eqref{lem15 1}, it can lead to
\begin{equation}
    \begin{aligned}
        & \left\|s_{k+1}-\mathcal{K} s_{k+1}\right\| \\
        & < \rho\left\|s_k-\mathcal{K} s_k\right\|+  L_2(L_3+1)\left\|y_{k+1}-y_k\right\|+2L_2\left\|u_k-\mathcal{K} u_k\right\| \\
        & \leq \rho\left\|s_k-\mathcal{K} s_k\right\|+ L_2(L_3+1)(\gamma+1)\left\|x_{k+1}-x_k\right\|\\
        & \quad+ L_2(L_3+1)\gamma\left\|x_{k}-x_{k-1}\right\|+2L_2\left\|u_k-\mathcal{K} u_k\right\|.
    \end{aligned}
    \label{lem15 3}
\end{equation}
Then by substituting \eqref{lem14} into \eqref{lem15 3}, we can obtain Lemma $15$.

\bibliographystyle{unsrt}
\bibliography{reference}  
\end{document}